\documentclass[12pt,english]{article}
\makeatletter

\usepackage{ifpdf}

\newif\ifpdf
\ifx\pdfoutput\undefined
  \pdffalse
\else
  \pdfoutput=1
  \pdftrue
\fi

\RequirePackage{xspace} %
\RequirePackage{subfigure} %

\ifpdf
  \RequirePackage[ pdftex, plainpages = false, pdfpagelabels,
                 pdfpagelayout = useoutlines,
                 bookmarks,
                 breaklinks = true,
                 linktocpage,
                 pagebackref,                      
                 colorlinks = true,
                 linkcolor = blue,
                 urlcolor  = blue,
                 citecolor = blue,
                 anchorcolor = blue,
                 hyperindex = true,
                 hyperfigures
                 ]{hyperref}
\usepackage{makeidx} 
\usepackage{setspace} 
\usepackage{rotating} 
\usepackage{ecltree}
\usepackage{epic}
\usepackage{supertabular}  
\usepackage{color}
\usepackage{exscale}
\usepackage{fontenc}
\usepackage{ifthen}
\usepackage{latexsym}
\usepackage{makeidx}
\usepackage{syntonly}
\usepackage{inputenc}
\usepackage{graphicx}
\usepackage{setspace}
\usepackage{caption2}
\usepackage[english]{babel}
\usepackage[square, comma,numbers,sort&compress]{natbib}
\usepackage{hypernat}
\usepackage{boxedminipage}
\usepackage{framed}
\usepackage{longtable}
\usepackage[all]{hypcap}
\usepackage{amsmath}
\usepackage{amssymb}
\usepackage{stmaryrd}
\usepackage{setspace}
\usepackage[bottom]{footmisc}

\newcommand{\note}[1]{\marginpar[left]{\singlespace \tiny #1}}

\renewcommand{\sectionmark}[1]%
      {\markright{\thesection\ #1}} 

\renewcommand{\note}[1]{}
\renewcommand\@biblabel[1]{$[{#1}]$}

\begin{document}

\title{Introduction to Tensor Calculus\\
\vspace{6cm}}

\author{Taha Sochi\thanks{Department of Physics \& Astronomy, University College London, Gower
Street, London, WC1E 6BT. Email: t.sochi@ucl.ac.uk.}\vspace{8cm}}

\maketitle
\pagebreak{}

\phantomsection \addcontentsline{toc}{section}{Preface}

\section*{Preface}

These are general notes on tensor calculus originated from a collection
of personal notes which I prepared some time ago for my own use and
reference when I was studying the subject. I decided to put them in
the public domain hoping they may be beneficial to some students in
their effort to learn this subject. Most of these notes were prepared
in the form of bullet points like tutorials and presentations and
hence some of them may be more concise than they should be. Moreover,
some notes may not be sufficiently thorough or general. However this
should be understandable considering the level and original purpose
of these notes and the desire for conciseness. There may also be some
minor repetition in some places for the purpose of gathering similar
items together, or emphasizing key points, or having self-contained
sections and units.

These notes, in my view, can be used as a short reference for an introductory
course on tensor algebra and calculus. I assume a basic knowledge
of calculus and linear algebra with some commonly used mathematical
terminology. I tried to be as clear as possible and to highlight the
key issues of the subject at an introductory level in a concise form.
I hope I have achieved some success in reaching these objectives at
least for some of my target audience. The present text is supposed
to be the first part of a series of documents about tensor calculus
for gradually increasing levels or tiers. I hope I will be able to
finalize and publicize the document for the next level in the near
future.

\pagebreak{}

\phantomsection \addcontentsline{toc}{section}{Contents}

\tableofcontents{}

\pagebreak{}

\section{Notation, Nomenclature and Conventions}

$\bullet$ In the present notes we largely follow certain conventions
and general notations; most of which are commonly used in the mathematical
literature although they may not be universally adopted. In the following
bullet points we outline these conventions and notations. We also
give initial definitions of the most basic terms and concepts in tensor
calculus; more thorough technical definitions will follow, if needed,
in the forthcoming sections.

$\bullet$ Scalars are algebraic objects which are uniquely identified
by their magnitude (absolute value) and sign ($\pm$), while vectors
are broadly geometric objects which are uniquely identified by their
magnitude (length) and direction in a presumed underlying space.

$\bullet$ At this early stage in these notes, we generically define
``tensor'' as an organized array of mathematical objects such as
numbers or functions.

$\bullet$ In generic terms, the rank of a tensor signifies the complexity
of its structure. Rank-0 tensors are called scalars while rank-1 tensors
are called vectors. Rank-2 tensors may be called dyads although this,
in common use, may be restricted to the outer product of two vectors
and hence is a special case of rank-2 tensors assuming it meets the
requirements of a tensor and hence transforms as a tensor. Like rank-2
tensors, rank-3 tensors may be called triads. Similar labels, which
are much less common in use, may be attached to higher rank tensors;
however, none will be used in the present notes. More generic names
for higher rank tensors, such as polyad, are also in use.

$\bullet$ In these notes we may use ``tensor'' to mean tensors
of all ranks including scalars (rank-0) and vectors (rank-1). We may
also use it as opposite to scalar and vector (i.e. tensor of rank-$n$
where $n>1$). In almost all cases, the meaning should be obvious
from the context.

$\bullet$ Non-indexed lower case light face Latin letters (e.g. $f$
and $h$) are used for scalars.

$\bullet$ Non-indexed (lower or upper case) bold face Latin letters
(e.g. $\mathbf{a}$ and $\mathbf{A}$) are used for vectors. The exception
to this is the basis vectors where indexed bold face lower or upper
case symbols are used. However, there should be no confusion or ambiguity
about the meaning of any one of these symbols.

$\bullet$ Non-indexed upper case bold face Latin letters (e.g. $\mathbf{A}$
and $\mathbf{B}$) are used for tensors (i.e. of rank $>1$).

$\bullet$ Indexed light face italic symbols (e.g. $a_{i}$ and $B_{i}^{jk}$)
are used to denote tensors of rank $>0$ in their explicit tensor
form (index notation). Such symbols may also be used to denote the
components of these tensors. The meaning is usually transparent and
can be identified from the context if not explicitly declared.

$\bullet$ Tensor indices in this document are lower case Latin letters
usually taken from the middle of the Latin alphabet like ($i,j,k$).
We also use numbered indices like ($i_{1},i_{2},\ldots,i_{k}$) when
the number of tensor indices is variable.

$\bullet$ The present notes are largely based on assuming an underlying
orthonormal Cartesian coordinate system. However, parts of which are
based on more general coordinate systems; in these cases this is stated
explicitly or made clear by the content and context.

$\bullet$ Mathematical identities and definitions may be denoted
by using the symbol `$\equiv$'. However, for simplicity we will use
in the present notes the equality sign ``='' to mark identities
and mathematical definitions as well as normal equalities.

$\bullet$ We use 2D, 3D and $n$D for  two-, three- and $n$-dimensional
spaces. We also use Eq./Eqs. to abbreviate Equation/Equations.

$\bullet$ Vertical bars are used to symbolize determinants while
square brackets are used for matrices.

$\bullet$ All tensors in the present notes are assumed to be real
quantities (i.e. have real rather than complex components).

$\bullet$ Partial derivative symbol with a subscript index (e.g.
$i$) is frequently used to denote the $i^{th}$ component of the
Cartesian gradient operator $\nabla$:
\begin{equation}
\partial_{i}=\nabla_{i}=\frac{\partial}{\partial x_{i}}
\end{equation}

$\bullet$ A comma preceding a subscript index (e.g. $,i$) is also
used to denote partial differentiation with respect to the $i^{th}$
spatial coordinate in Cartesian systems, e.g.
\begin{equation}
A_{,i}=\frac{\partial A}{\partial x_{i}}
\end{equation}

$\bullet$ Partial derivative symbol with a spatial subscript, rather
than an index, are used to denote partial differentiation with respect
to that spatial variable. For instance
\begin{equation}
\partial_{r}=\nabla_{r}=\frac{\partial}{\partial r}
\end{equation}
is used for the partial derivative with respect to the radial coordinate
in spherical coordinate systems identified by ($r,\theta,\phi$) spatial
variables.

$\bullet$ Partial derivative symbol with repeated double index is
used to denote the Laplacian operator:
\begin{equation}
\partial_{ii}=\partial_{i}\partial_{i}=\nabla^{2}=\Delta\label{eqLaplacianSymbol}
\end{equation}
The notation is not affected by using repeated double index other
than $i$ (e.g. $\partial_{jj}$ or $\partial_{kk}$). The following
notations:
\begin{equation}
\partial_{ii}^{2}\,\,\,\,\,\,\,\,\,\,\,\,\,\,\,\partial^{2}\,\,\,\,\,\,\,\,\,\,\,\,\,\,\,\partial_{i}\partial^{i}
\end{equation}
are also used in the literature of tensor calculus to symbolize the
Laplacian operator. However, these notations will not be used in the
present notes.

$\bullet$ We follow the common convention of using a subscript semicolon
preceding a subscript index (e.g. $A_{kl;i}$) to symbolize covariant
differentiation with respect to the $i^{th}$ coordinate (see $\S$
\ref{secCovariantDifferentiation}). The semicolon notation may also
be attached to the normal differential operators to indicate covariant
differentiation (e.g. $\nabla_{;i}$ or $\partial_{;i}$ to indicate
covariant differentiation with respect to the index $i$).

$\bullet$ All transformation equations in these notes are assumed
continuous and real, and all derivatives are continuous in their domain
of variables.

$\bullet$ Based on the continuity condition of the differentiable
quantities, the individual differential operators in the mixed partial
derivatives are commutative, that is
\begin{equation}
\partial_{i}\partial_{j}=\partial_{j}\partial_{i}
\end{equation}

$\bullet$ A permutation of a set of objects, which are normally numbers
like $\left(1,2,\ldots,n\right)$ or symbols like $\left(i,j,k\right)$,
is a particular ordering or arrangement of these objects. An even
permutation is a permutation resulting from an even number of single-step
exchanges (also known as transpositions) of neighboring objects starting
from a presumed original permutation of these objects. Similarly,
an odd permutation is a permutation resulting from an odd number of
such exchanges. It has been shown that when a transformation from
one permutation to another can be done in different ways, possibly
with different numbers of exchanges, the parity of all these possible
transformations is the same, i.e. all even or all odd, and hence there
is no ambiguity in characterizing the transformation from one permutation
to another by the parity alone.

$\bullet$ We normally use indexed square brackets (e.g. $\left[\mathbf{A}\right]_{i}$
and $\left[\mathbf{\nabla}f\right]_{i}$) to denote the $i^{th}$
component of vectors, tensors and operators in their symbolic or vector
notation.

$\bullet$ In general terms, a transformation from an $n$D space
to another $n$D space is a correlation that maps a point from the
first space (original) to a point in the second space (transformed)
where each point in the original and transformed spaces is identified
by $n$ independent variables or coordinates. To distinguish between
the two sets of coordinates in the two spaces, the coordinates of
the points in the transformed space may be notated with barred symbols,
e.g. ($\bar{x}^{1},\bar{x}^{2},\ldots,\bar{x}^{n}$) or ($\bar{x}_{1},\bar{x}_{2},\ldots,\bar{x}_{n}$)
where the superscripts and subscripts are indices, while the coordinates
of the points in the original space are notated with unbarred symbols,
e.g. ($x^{1},x^{2},\ldots,x^{n}$) or ($x_{1},x_{2},\ldots,x_{n}$).
Under certain conditions, such a transformation is unique and hence
an inverse transformation from the transformed to the original space
is also defined. Mathematically, each one of the direct and inverse
transformation can be regarded as a mathematical correlation expressed
by a set of equations in which each coordinate in one space is considered
as a function of the coordinates in the other space. Hence the transformations
between the two sets of coordinates in the two spaces can by expressed
mathematically by the following two sets of independent relations:
\begin{equation}
\bar{x}^{i}=\bar{x}^{i}(x^{1},x^{2},\ldots,x^{n})\,\,\,\,\,\,\,\,\,\,\&\,\,\,\,\,\,\,\,\,\,x^{i}=x^{i}(\bar{x}^{1},\bar{x}^{2},\ldots,\bar{x}^{n})
\end{equation}
where $i=1,2,\ldots,n$. An alternative to viewing the transformation
as a mapping between two different spaces is to view it as being correlating
the same point in the same space but observed from two different coordinate
frames of reference which are subject to a similar transformation.

$\bullet$ Coordinate transformations are described as ``proper''
when they preserve the handedness (right- or left-handed) of the coordinate
system and ``improper'' when they reverse the handedness. Improper
transformations involve an odd number of coordinate axes inversions
through the origin.

$\bullet$ Inversion of axes may be called improper rotation while
ordinary rotation is described as proper rotation.

$\bullet$ Transformations can be active, when they change the state
of the observed object (e.g. translating the object in space), or
passive when they are based on keeping the state of the object and
changing the state of the coordinate system from which the object
is observed. Such distinction is based on an implicit assumption of
a more general frame of reference in the background.

$\bullet$ Finally, tensor calculus is riddled with conflicting conventions
and terminology. In this text we will try to use what we believe to
be the most common, clear or useful of all of these.

\pagebreak{}

\section{Preliminaries}

\subsection{Introduction}

$\bullet$ A tensor is an array of mathematical objects (usually numbers
or functions) which transforms according to certain rules under coordinates
change. In a $d$-dimensional space, a tensor of rank-$n$ has $d^{n}$
components which may be specified with reference to a given coordinate
system. Accordingly, a scalar, such as temperature, is a rank-0 tensor
with (assuming 3D space) $3^{0}=1$ component, a vector, such as force,
is a rank-1 tensor with $3^{1}=3$ components, and stress is a rank-2
tensor with $3^{2}=9$ components.

$\bullet$ The term ``tensor'' was originally derived from the Latin
word ``tensus'' which means tension or stress since one of the first
uses of tensors was related to the mathematical description of mechanical
stress.

$\bullet$ The $d^{n}$ components of a tensor are identified by $n$
distinct integer indices (e.g. $i,j,k$) which are attached, according
to the commonly-employed tensor notation, as superscripts or subscripts
or a mix of these to the right side of the symbol utilized to label
the tensor (e.g. $A_{ijk}$, $A^{ijk}$ and $A_{i}^{jk}$). Each tensor
index takes all the values over a predefined range of dimensions such
as 1 to $d$ in the above example of a $d$-dimensional space. In
general, all tensor indices have the same range, i.e. they are uniformly
dimensioned.

$\bullet$ When the range of tensor indices is not stated explicitly,
it is usually assumed to have the values ($1,2,3$). However, the
range must be stated explicitly or implicitly to avoid ambiguity.

$\bullet$ The characteristic property of tensors is that they satisfy
the principle of invariance under certain coordinate transformations.
Therefore, formulating the fundamental physical laws in a tensor form
ensures that they are form-invariant; hence they are objectively-representing
the physical reality and do not depend on the observer. Having the
same form in different coordinate systems may be labeled as being
``covariant'' but this word is also used for a different meaning
in tensor calculus as explained in $\S$ \ref{subCovariantContravariant}.

$\bullet$ ``Tensor term'' is a product of tensors including scalars
and vectors.

$\bullet$ ``Tensor expression'' is an algebraic sum (or more generally
a linear combination) of tensor terms which may be a trivial sum in
the case of a single term.

$\bullet$ ``Tensor equality'' (symbolized by `$=$') is an equality
of two tensor terms and/or expressions. A special case of this is
tensor identity which is an equality of general validity (the symbol
`$\equiv$' may be used for identity as well as for definition).

$\bullet$ The order of a tensor is identified by the number of its
indices (e.g. $A_{jk}^{i}$ is a tensor of order 3) which normally
identifies the tensor rank as well. However, when contraction (see
$\S$ \ref{subsubContraction}) takes place once or more, the order
of the tensor is not affected but its rank is reduced by two for each
contraction operation.\footnote{In the literature of tensor calculus, rank and order of tensors are
generally used interchangeably; however some authors differentiate
between the two as they assign order to the total number of indices,
including repetitive indices, while they keep rank to the number of
free indices. We think the latter is better and hence we follow this
convention in the present text. }

$\bullet$ ``Zero tensor'' is a tensor whose all components are
zero.

$\bullet$ ``Unit tensor'' or ``unity tensor'', which is usually
defined for rank-2 tensors, is a tensor whose all elements are zero
except the ones with identical values of all indices which are assigned
the value 1.

$\bullet$ While tensors of rank-0 are generally represented in a
common form of light face non-indexed symbols, tensors of rank $\ge1$
are represented in several forms and notations, the main ones are
the index-free notation, which may also be called direct or symbolic
or Gibbs notation, and the indicial notation which is also called
index or component or tensor notation. The first is a geometrically
oriented notation with no reference to a particular reference frame
and hence it is intrinsically invariant to the choice of coordinate
systems, whereas the second takes an algebraic form based on components
identified by indices and hence the notation is suggestive of an underlying
coordinate system, although being a tensor makes it form-invariant
under certain coordinate transformations and therefore it possesses
certain invariant properties. The index-free notation is usually identified
by using bold face symbols, like $\mathbf{a}$ and $\mathbf{B}$,
while the indicial notation is identified by using light face indexed
symbols such as $a^{i}$ and $B_{ij}$.

\subsection{General Rules\label{subGeneralRules} }

$\bullet$ An index that occurs once in a tensor term is a ``free
index''.

$\bullet$ An index that occurs twice in a tensor term is a ``dummy''
or ``bound'' index.

$\bullet$ No index is allowed to occur more than twice in a legitimate
tensor term.\footnote{We adopt this assertion, which is common in the literature of tensor
calculus, as we think it is suitable for this level. However, there
are many instances in the literature of tensor calculus where indices
are repeated more than twice in a single term. The bottom line is
that as long as the tensor expression makes sense and the intention
is clear, such repetitions should be allowed with no need in our view
to take special precaution like using parentheses. In particular,
the summation convention will not apply automatically in such cases
although summation on such indices can be carried out explicitly,
by using the summation symbol $\sum$, or by special declaration of
such intention similar to the summation convention. Anyway, in the
present text we will not use indices repeated more than twice in a
single term.}

$\bullet$ A free index should be understood to vary over its range
(e.g. $1,\ldots,n$) and hence it can be interpreted as saying ``for
all components represented by the index''. Therefore a free index
represents a number of terms or expressions or equalities equal to
the number of allowed values of its range. For example, when $i$
and $j$ can vary over the range $1,\ldots,n$ the following expression
\begin{equation}
A_{i}+B_{i}
\end{equation}
represents $n$ separate expressions while the following equation
\begin{equation}
A_{i}^{j}=B_{i}^{j}
\end{equation}
represents $n\times n$ separate equations.

$\bullet$ According to the ``summation convention'', which is widely
used in the literature of tensor calculus including in the present
notes, dummy indices imply summation over their range, e.g. for an
$n$D space
\begin{equation}
A^{i}B_{i}\equiv\sum_{i=1}^{n}A^{i}B_{i}=A^{1}B_{1}+A^{2}B_{2}+\ldots+A^{n}B_{n}
\end{equation}
\begin{equation}
\delta_{ij}A^{ij}\equiv\sum_{i=1}^{n}\sum_{j=1}^{n}\delta_{ij}A^{ij}
\end{equation}
\begin{equation}
\epsilon_{ijk}A^{ij}B^{k}\equiv\sum_{i=1}^{n}\sum_{j=1}^{n}\sum_{k=1}^{n}\epsilon_{ijk}A^{ij}B^{k}
\end{equation}

$\bullet$ When dummy indices do not imply summation, the situation
must be clarified by enclosing such indices in parentheses or by underscoring
or by using upper case letters (with declaration of these conventions)
or by adding a clarifying comment like ``no summation on repeated
indices''.

$\bullet$ Tensors with subscript indices, like $A_{ij}$, are called
covariant, while tensors with superscript indices, like $A^{k}$,
are called contravariant. Tensors with both types of indices, like
$A_{lk}^{lmn}$, are called mixed type. More details about this will
follow in $\S$ \ref{subCovariantContravariant}.

$\bullet$ Subscript indices, rather than subscripted tensors, are
also dubbed ``covariant'' and superscript indices are dubbed ``contravariant''.

$\bullet$ Each tensor index should conform to one of the variance
transformation rules as given by Eqs. \ref{eqCovariant} and \ref{eqContravariant},
i.e. it is either covariant or contravariant.

$\bullet$ For orthonormal Cartesian coordinate systems, the two variance
types (i.e. covariant and contravariant) do not differ because the
metric tensor is given by the Kronecker delta (refer to $\S$ \ref{secMetricTensor}
and \ref{subKronecker}) and hence any index can be upper or lower
although it is common to use lower indices in such cases.

$\bullet$ For tensor invariance, a pair of dummy indices should in
general be complementary in their variance type, i.e. one covariant
and the other contravariant. However, for orthonormal Cartesian systems
the two are the same and hence when both dummy indices are covariant
or both are contravariant it should be understood as an indication
that the underlying coordinate system is orthonormal Cartesian if
the possibility of an error is excluded.

$\bullet$ As indicated earlier, tensor order is equal to the number
of its indices while tensor rank is equal to the number of its free
indices; hence vectors (terms, expressions and equalities) are represented
by a single free index and rank-2 tensors are represented by two free
indices. The dimension of a tensor is determined by the range taken
by its indices.

$\bullet$ The rank of all terms in legitimate tensor expressions
and equalities must be the same.

$\bullet$ Each term in valid tensor expressions and equalities must
have the same set of free indices (e.g. $i,j,k$).

$\bullet$ A free index should keep its variance type in every term
in valid tensor expressions and equations, i.e. it must be covariant
in all terms or contravariant in all terms.

$\bullet$ While free indices should be named uniformly in all terms
of tensor expressions and equalities, dummy indices can be named in
each term independently, e.g.
\begin{equation}
A_{ik}^{i}+B_{jk}^{j}+C_{lmk}^{lm}
\end{equation}

$\bullet$ A free index in an expression or equality can be renamed
uniformly using a different symbol, as long as this symbol is not
already in use, assuming that both symbols vary over the same range,
i.e. have the same dimension.

$\bullet$ Examples of legitimate tensor terms, expressions and equalities:
\begin{equation}
A_{ij}^{ij},\,\,\,\,\,\,\,\,\,\,A_{m}^{im}+B_{nk}^{ink},\,\,\,\,\,\,\,\,\,\,C_{ij}=A_{ij}-B_{ij},\,\,\,\,\,\,\,\,\,\,a=B_{j}^{j}
\end{equation}

$\bullet$ Examples of illegitimate tensor terms, expressions and
equalities:
\begin{equation}
B_{i}^{ii},\,\,\,\,\,\,\,\,\,\,A_{i}+B_{ij},\,\,\,\,\,\,\,\,\,\,A^{i}+B^{j},\,\,\,\,\,\,\,\,\,\,A_{i}-B^{i},\,\,\,\,\,\,\,\,\,\,A_{i}^{i}=B_{i},
\end{equation}

$\bullet$ Indexing is generally distributive over the terms of tensor
expressions and equalities, e.g.
\begin{equation}
\left[\mathbf{A}+\mathbf{B}\right]_{i}=\left[\mathbf{A}\right]_{i}+\left[\mathbf{B}\right]_{i}\label{eqIndexDistributive1}
\end{equation}
and
\begin{equation}
\left[\mathbf{A}=\mathbf{B}\right]_{i}\,\,\,\,\,\,\,\,\,\,\Longleftrightarrow\,\,\,\,\,\,\,\,\,\,\left[\mathbf{A}\right]_{i}=\left[\mathbf{B}\right]_{i}
\end{equation}

$\bullet$ Unlike scalars and tensor components, which are essentially
scalars in a generic sense, operators cannot in general be freely
reordered in tensor terms, therefore
\begin{equation}
fh=hf\,\,\,\,\,\,\,\,\,\,\&\,\,\,\,\,\,\,\,\,\,A_{i}B^{i}=B^{i}A_{i}
\end{equation}
but
\begin{equation}
\partial_{i}A_{i}\ne A_{i}\partial_{i}
\end{equation}

$\bullet$ Almost all the identities in the present notes which are
given in a covariant or a contravariant or a mixed form are similarly
valid for the other forms unless it is stated otherwise. The objective
of reporting in only one form is conciseness and to avoid unnecessary
repetition.

\subsection{Examples of Tensors of Different Ranks}

$\bullet$ Examples of rank-0 tensors (scalars) are energy, mass,
temperature, volume and density. These are totally identified by a
single number regardless of any coordinate system and hence they are
invariant under coordinate transformations.

$\bullet$ Examples of rank-1 tensors (vectors) are displacement,
force, electric field, velocity and acceleration. These need for their
complete identification a number, representing their magnitude, and
a direction representing their geometric orientation within their
space. Alternatively, they can be uniquely identified by a set of
numbers, equal to the number of dimensions of the underlying space,
in reference to a particular coordinate system and hence this identification
is system-dependent although they still have system-invariant properties
such as length.

$\bullet$ Examples of rank-2 tensors are Kronecker delta (see $\S$
\ref{subKronecker}), stress, strain, rate of strain and inertia tensors.
These require for their full identification a set of numbers each
of which is associated with two directions.

$\bullet$ Examples of rank-3 tensors are the Levi-Civita tensor (see
$\S$ \ref{subPermutation}) and the tensor of piezoelectric moduli.

$\bullet$ Examples of rank-4 tensors are the elasticity or stiffness
tensor, the compliance tensor and the fourth-order moment of inertia
tensor.

$\bullet$ Tensors of high ranks are relatively rare in science and
engineering.

$\bullet$ Although rank-0 and rank-1 tensors are, respectively, scalars
and vectors, not all scalars and vectors (in their generic sense)
are tensors of these ranks. Similarly, rank-2 tensors are normally
represented by matrices but not all matrices represent tensors.

\subsection{Applications of Tensors}

$\bullet$ Tensor calculus is very powerful mathematical tool. Tensor
notation and techniques are used in many branches of science and engineering
such as fluid mechanics, continuum mechanics, general relativity and
structural engineering. Tensor calculus is used for elegant and compact
formulation and presentation of equations and identities in mathematics,
science and engineering. It is also used for algebraic manipulation
of mathematical expressions and proving identities in a neat and succinct
way (refer to $\S$ \ref{subProvingIdentities}).

$\bullet$ As indicated earlier, the invariance of tensor forms serves
a theoretically and practically important role by allowing the formulation
of physical laws in coordinate-free forms.

\subsection{Types of Tensors}

$\bullet$ In the following subsections we introduce a number of tensor
types and categories and highlight their main characteristics and
differences. These types and categories are not mutually exclusive
and hence they overlap in general; moreover they may not be exhaustive
in their classes as some tensors may not instantiate any one of a
complementary set of types such as being symmetric or anti-symmetric.

\subsubsection{Covariant and Contravariant Tensors\label{subCovariantContravariant}}

$\bullet$ These are the main types of tensor with regard to the rules
of their transformation between different coordinate systems.

$\bullet$ Covariant tensors are notated with subscript indices (e.g.
$A_{i}$) while contravariant tensors are notated with superscript
indices (e.g. $A^{ij}$).

$\bullet$ A covariant tensor is transformed according to the following
rule
\begin{equation}
\bar{A}{}_{i}=\frac{\partial x^{j}}{\partial\bar{x}{}^{i}}A_{j}\label{eqCovariant}
\end{equation}
while a contravariant tensor is transformed according to the following
rule
\begin{equation}
\bar{A}{}^{i}=\frac{\partial\bar{x}{}^{i}}{\partial x{}^{j}}A^{j}\label{eqContravariant}
\end{equation}
where the barred and unbarred symbols represent the same mathematical
object (tensor or coordinate) in the transformed and original coordinate
systems respectively.

$\bullet$ An example of covariant tensors is the gradient of a scalar
field.

$\bullet$ An example of contravariant tensors is the displacement
vector.

$\bullet$ Some tensors have mixed variance type, i.e. they are covariant
in some indices and contravariant in others. In this case the covariant
variables are indexed with subscripts while the contravariant variables
are indexed with superscripts, e.g. $A_{i}^{j}$ which is covariant
in $i$ and contravariant in $j$.

$\bullet$ A mixed type tensor transforms covariantly in its covariant
indices and contravariantly in its contravariant indices, e.g.
\begin{equation}
\bar{A}{}_{\,\,m}^{l\,\,\,\,n}=\frac{\partial\bar{x}{}^{l}}{\partial x^{i}}\frac{\partial x^{j}}{\partial\bar{x}{}^{m}}\frac{\partial\bar{x}{}^{n}}{\partial x^{k}}A{}_{\,\,\,j}^{i\,\,\,k}
\end{equation}

$\bullet$ To clarify the pattern of mathematical transformation of
tensors, we explain step-by-step the practical rules to follow in
writing tensor transformation equations between two coordinate systems,
unbarred and barred, where for clarity we color the symbols of the
tensor and the coordinates belonging to the unbarred system with blue
while we use red to mark the symbols belonging to the barred system.
Since there are three types of tensors: covariant, contravariant and
mixed, we use three equations in each step. In this demonstration
we use rank-4 tensors as examples since this is sufficiently general
and hence adequate to elucidate the rules for tensors of any rank.
The demonstration is based on the assumption that the transformation
is taking place from the unbarred system to the barred system; the
same rules should apply for the opposite transformation from the barred
system to the unbarred system. We use the sign `$\circeq$' for the
equality in the transitional steps to indicate that the equalities
are under construction and are not complete.

We start by the very generic equations between the barred tensor ${\color{red}\bar{A}}$
and the unbarred tensor ${\color{blue}A}$ for the three types:
\begin{eqnarray}
{\color{blue}{\color{red}\bar{A}}} & \circeq & {\color{blue}A}\,\,\,\,\,\,\,\,\,\,\,\,\,\,\,\,\,\,\,\,\text{(covariant)}\nonumber \\
{\color{blue}{\color{red}\bar{A}}} & \circeq & {\color{blue}A}\,\,\,\,\,\,\,\,\,\,\,\,\,\,\,\,\,\,\,\,\text{(contravariant)}\\
{\color{blue}{\color{red}\bar{A}}} & \circeq & {\color{blue}A}\,\,\,\,\,\,\,\,\,\,\,\,\,\,\,\,\,\,\,\,\text{(mixed)}\nonumber
\end{eqnarray}
We assume that the barred tensor and its coordinates are indexed with
$ijkl$ and the unbarred are indexed with $npqr$, so we add these
indices in their presumed order and position (lower or upper) paying
particular attention to the order in the mixed type:
\begin{eqnarray}
{\color{red}\bar{A}_{ijkl}} & \circeq & {\color{blue}A_{npqr}}\nonumber \\
{\color{red}\bar{A}^{ijkl}} & \circeq & {\color{blue}A^{npqr}}\\
{\color{red}\bar{A}_{\,\,\,kl}^{ij}} & \circeq & {\color{blue}A_{\,\,\,\,\,qr}^{np}}\nonumber
\end{eqnarray}
Since the barred and unbarred tensors are of the same type, as they
represent the same tensor in two coordinate systems,\footnote{Similar basis vectors are assumed.}
the indices on the two sides of the equalities should match in their
position and order. We then insert a number of partial differential
operators on the right hand side of the equations equal to the rank
of these tensors, which is 4 in our example. These operators represent
the transformation rules for each pair of corresponding coordinates
one from the barred and one from the unbarred:
\begin{eqnarray}
{\color{red}\bar{A}_{ijkl}} & \circeq\frac{\partial}{\partial}\,\frac{\partial}{\partial}\,\frac{\partial}{\partial}\,\frac{\partial}{\partial} & {\color{blue}A_{npqr}}\nonumber \\
{\color{red}\bar{A}^{ijkl}} & \circeq\frac{\partial}{\partial}\,\frac{\partial}{\partial}\,\frac{\partial}{\partial}\,\frac{\partial}{\partial} & {\color{blue}A^{npqr}}\\
{\color{red}\bar{A}_{\,\,\,kl}^{ij}} & \circeq\frac{\partial}{\partial}\,\frac{\partial}{\partial}\,\frac{\partial}{\partial}\,\frac{\partial}{\partial} & {\color{blue}A_{\,\,\,\,\,qr}^{np}}\nonumber
\end{eqnarray}
Now we insert the coordinates of the barred system into the partial
differential operators noting that (i) the positions of any index
on the two sides should match, i.e. both upper or both lower, since
they are free indices in different terms of tensor equalities, (ii)
a superscript index in the denominator of a partial derivative is
in lieu of a covariant index in the numerator\footnote{The use of upper indices in the denominator of partial derivatives,
which is common in this type of equations, is to indicate the fact
that the coordinates and their differentials transform contravariantly.}, and (iii) the order of the coordinates should match the order of
the indices in the tensor:
\begin{eqnarray}
{\color{red}\bar{A}_{ijkl}} & \circeq\frac{\partial}{\partial{\color{red}x^{i}}}\,\frac{\partial}{\partial{\color{red}x^{j}}}\,\frac{\partial}{\partial{\color{red}x^{k}}}\,\frac{\partial}{\partial{\color{red}x^{l}}} & {\color{blue}A_{npqr}}\nonumber \\
{\color{red}\bar{A}^{ijkl}} & \circeq\frac{\partial{\color{red}x^{i}}}{\partial}\,\frac{\partial{\color{red}x^{j}}}{\partial}\,\frac{\partial{\color{red}x^{k}}}{\partial}\,\frac{\partial{\color{red}x^{l}}}{\partial} & {\color{blue}A^{npqr}}\\
{\color{red}\bar{A}_{\,\,\,kl}^{ij}} & \circeq\frac{\partial{\color{red}x^{i}}}{\partial}\,\frac{\partial{\color{red}x^{j}}}{\partial}\,\frac{\partial}{\partial{\color{red}x^{k}}}\,\frac{\partial}{\partial{\color{red}x^{l}}} & {\color{blue}A_{\,\,\,\,\,qr}^{np}}\nonumber
\end{eqnarray}
For consistency, these coordinates should be barred as they belong
to the barred tensor; hence we add bars:
\begin{eqnarray}
{\color{red}\bar{A}_{ijkl}} & \circeq\frac{\partial}{\partial{\color{red}\bar{x}^{i}}}\,\frac{\partial}{\partial{\color{red}\bar{x}^{j}}}\,\frac{\partial}{\partial{\color{red}\bar{x}^{k}}}\,\frac{\partial}{\partial{\color{red}\bar{x}^{l}}} & {\color{blue}A_{npqr}}\nonumber \\
{\color{red}\bar{A}^{ijkl}} & \circeq\frac{\partial{\color{red}\bar{x}^{i}}}{\partial}\,\frac{\partial{\color{red}\bar{x}^{j}}}{\partial}\,\frac{\partial{\color{red}\bar{x}^{k}}}{\partial}\,\frac{\partial{\color{red}\bar{x}^{l}}}{\partial} & {\color{blue}A^{npqr}}\\
{\color{red}\bar{A}_{\,\,\,kl}^{ij}} & \circeq\frac{\partial{\color{red}\bar{x}^{i}}}{\partial}\,\frac{\partial{\color{red}\bar{x}^{j}}}{\partial}\,\frac{\partial}{\partial{\color{red}\bar{x}^{k}}}\,\frac{\partial}{\partial{\color{red}\bar{x}^{l}}} & {\color{blue}A_{\,\,\,\,\,qr}^{np}}\nonumber
\end{eqnarray}
Finally, we insert the coordinates of the unbarred system into the
partial differential operators noting that (i) the positions of the
repeated indices on the same side should be opposite, i.e. one upper
and one lower, since they are dummy indices and hence the position
of the index of the unbarred coordinate should be opposite to its
position in the unbarred tensor, (ii) an upper index in the denominator
is in lieu of a lower index in the numerator, and (iii) the order
of the coordinates should match the order of the indices in the tensor:
\begin{eqnarray}
{\color{red}\bar{A}_{ijkl}} & =\frac{\partial{\color{blue}x^{n}}}{\partial{\color{red}\bar{x}^{i}}}\,\frac{\partial{\color{blue}x^{p}}}{\partial{\color{red}\bar{x}^{j}}}\,\frac{\partial{\color{blue}x^{q}}}{\partial{\color{red}\bar{x}^{k}}}\,\frac{\partial{\color{blue}x^{r}}}{\partial{\color{red}\bar{x}^{l}}} & {\color{blue}A_{npqr}}\nonumber \\
{\color{red}\bar{A}^{ijkl}} & =\frac{\partial{\color{red}\bar{x}^{i}}}{\partial{\color{blue}x^{n}}}\,\frac{\partial{\color{red}\bar{x}^{j}}}{\partial{\color{blue}x^{p}}}\,\frac{\partial{\color{red}\bar{x}^{k}}}{\partial{\color{blue}x^{q}}}\,\frac{\partial{\color{red}\bar{x}^{l}}}{\partial{\color{blue}x^{r}}} & {\color{blue}A^{npqr}}\\
{\color{red}\bar{A}_{\,\,\,kl}^{ij}} & =\frac{\partial{\color{red}\bar{x}^{i}}}{\partial{\color{blue}x^{n}}}\,\frac{\partial{\color{red}\bar{x}^{j}}}{\partial{\color{blue}x^{p}}}\,\frac{\partial{\color{blue}x^{q}}}{\partial{\color{red}\bar{x}^{k}}}\,\frac{\partial{\color{blue}x^{r}}}{\partial{\color{red}\bar{x}^{l}}} & {\color{blue}A_{\,\,\,\,\,qr}^{np}}\nonumber
\end{eqnarray}
We also replaced the `$\circeq$' sign in the final set of equations
with the strict equality sign `=' as the equations now are complete.

$\bullet$ A tensor of $m$ contravariant indices and $n$ covariant
indices may be called type ($m,n$) tensor, e.g. $A_{ij}^{k}$ is
a type ($1,2$) tensor. When one or both variance types are absent,
zero is used to refer to the absent type in this notation, e.g. $B^{ik}$
is a type ($2,0$) tensor.

$\bullet$ The covariant and contravariant types of a tensor are linked
through the metric tensor (refer to $\S$ \ref{secMetricTensor}).

$\bullet$ For orthonormal Cartesian systems there is no difference
between covariant and contravariant tensors, and hence the indices
can be upper or lower.

$\bullet$ The vectors providing the basis set (not necessarily of
unit length or mutually orthogonal) for a coordinate system are of
covariant type when they are tangent to the coordinate axes, and they
are of contravariant type when they are perpendicular to the local
surfaces of constant coordinates. These two sets are identical for
orthonormal Cartesian systems.

$\bullet$ Formally, the covariant and contravariant basis vectors
are given respectively by:
\begin{equation}
\mathbf{E}_{i}=\frac{\partial\mathbf{r}}{\partial u^{i}}\,\,\,\,\,\,\,\,\,\,\,\,\,\,\,\,\,\,\,\,\&\,\,\,\,\,\,\,\,\,\,\,\,\,\,\,\,\,\,\,\,\mathbf{E}^{i}=\nabla u^{i}
\end{equation}
where $\mathbf{r}$ is the position vector in Cartesian coordinates
and $u^{i}$ is a generalized curvilinear coordinate. As indicated
already, a superscript in the denominator of partial derivatives is
equivalent to a subscript in the numerator.

$\bullet$ In general, the covariant and contravariant basis vectors
are not mutually orthogonal or of unit length; however the two sets
are reciprocal systems and hence they satisfy the following reciprocity
relation:
\begin{equation}
\mathbf{E}_{i}\cdot\mathbf{E}^{j}=\delta_{i}^{j}
\end{equation}
where $\delta_{i}^{j}$ is the Kronecker delta (refer to $\S$ \ref{subKronecker}).

$\bullet$ A vector can be represented either by covariant components
with contravariant coordinate basis vectors or by contravariant components
with covariant coordinate basis vectors. For example, a vector $\mathbf{A}$
can be expressed as
\begin{equation}
\mathbf{A}=A_{i}\mathbf{E}^{i}\,\,\,\,\,\,\,\,\,\,\,\,\,\,\,\,\,\,\,\,\text{or}\,\,\,\,\,\,\,\,\,\,\,\,\,\,\,\,\,\,\,\,\mathbf{A}=A^{i}\mathbf{E}_{i}
\end{equation}
where $\mathbf{E}^{i}$ and $\mathbf{E}_{i}$ are the contravariant
and covariant basis vectors respectively. The use of the covariant
or contravariant form of the vector representation is a matter of
choice and convenience.

$\bullet$ More generally, a tensor of any rank ($\ge1$) can be represented
covariantly using contravariant basis tensors of that rank, or contravariantly
using covariant basis tensors, or in a mixed form using a mixed basis
of opposite type. For example, a rank-2 tensor $\mathbf{A}$ can be
written as:
\begin{equation}
\mathbf{A}=A_{ij}\mathbf{E}^{i}\mathbf{E}^{j}=A^{ij}\mathbf{E}_{i}\mathbf{E}_{j}=A_{i}^{\,\,j}\mathbf{E}^{i}\mathbf{E}_{j}
\end{equation}
where $\mathbf{E}^{i}\mathbf{E}^{j}$, $\mathbf{E}_{i}\mathbf{E}_{j}$
and $\mathbf{E}^{i}\mathbf{E}_{j}$ are dyadic products (refer to
$\S$ \ref{subTensorMultiplication}).

\subsubsection{True and Pseudo Tensors}

$\bullet$ These are also called polar and axial tensors respectively
although it is more common to use the latter terms for vectors. Pseudo
tensors may also be called tensor densities.

$\bullet$ True tensors are proper (or ordinary) tensors and hence
they are invariant under coordinate transformations, while pseudo
tensors are not proper tensors since they do not transform invariantly
as they acquire a minus sign under improper orthogonal transformations
which involve inversion of coordinate axes through the origin with
a change of system handedness.

$\bullet$ Because true and pseudo tensors have different mathematical
properties and represent different types of physical entities, the
terms of consistent tensor expressions and equations should be uniform
in their true and pseudo type, i.e. all terms are true or all are
pseudo.

$\bullet$ The direct product (refer to $\S$ \ref{subTensorMultiplication})
of even number of pseudo tensors is a proper tensor, while the direct
product of odd number of pseudo tensors is a pseudo tensor. The direct
product of true tensors is obviously a true tensor.

$\bullet$ The direct product of a mix of true and pseudo tensors
is a true or pseudo tensor depending on the number of pseudo tensors
involved in the product as being even or odd respectively.

$\bullet$ Similar rules to those of direct product apply to cross
products (including curl operations) involving tensors (usually of
rank-1) with the addition of a pseudo factor for each cross product
operation. This factor is contributed by the permutation tensor (see
$\S$ \ref{subPermutation}) which is implicit in the definition of
the cross product (see Eqs. \ref{EqCrossProduct} and \ref{EqCurl}).

$\bullet$ In summary, what determines the tensor type (true or pseudo)
of the tensor terms involving direct\footnote{Inner product (see $\S$ \ref{secInnerProduct}) is the result of
a direct product operation followed by a contraction (see $\S$ \ref{subsubContraction})
and hence it is a direct product in this context.} and cross products is the parity of the multiplicative factors of
pseudo type plus the number of cross product operations involved since
each cross product contributes an $\epsilon$ factor.

$\bullet$ Examples of true scalars are temperature, mass and the
dot product of two polar or two axial vectors, while examples of pseudo
scalars are the dot product of an axial vector and a polar vector
and the scalar triple product of polar vectors.

$\bullet$ Examples of polar vectors are displacement and acceleration,
while examples of axial vectors are angular velocity and cross product
of polar vectors in general (including curl operation on polar vectors)
due to the involvement of the permutation symbol which is a pseudo
tensor (refer to $\S$ \ref{subPermutation}). The essence of this
distinction is that the direction of a pseudo vector depends on the
observer choice of the handedness of the coordinate system whereas
the direction of a proper vector is independent of such choice.

$\bullet$ Examples of proper tensors of rank-2 are stress and rate
of strain tensors, while examples of pseudo tensors of rank-2 are
direct products of two vectors: one polar and one axial.

$\bullet$ Examples of proper tensors of higher ranks are piezoelectric
moduli tensor (rank-3) and elasticity tensor (rank-4), while examples
of pseudo tensors of higher ranks are the permutation tensor of these
ranks.

\subsubsection{Absolute and Relative Tensors\label{subAbsoluteRelative}}

$\bullet$ Considering an arbitrary transformation from a general
coordinate system to another, a relative tensor of weight $w$ is
defined by the following tensor transformation:
\begin{equation}
\bar{A}_{\,\,\,\,\,\,\,\,\,\,\,lm\ldots n}^{ij\ldots k}=\left|\frac{\partial x}{\partial\bar{x}}\right|^{w}\frac{\partial\bar{x}^{i}}{\partial x^{a}}\frac{\partial\bar{x}^{j}}{\partial x^{b}}\cdots\frac{\partial\bar{x}^{k}}{\partial x^{c}}\frac{\partial x^{d}}{\partial\bar{x}^{l}}\frac{\partial x^{e}}{\partial\bar{x}^{m}}\cdots\frac{\partial x^{f}}{\partial\bar{x}^{n}}A_{\,\,\,\,\,\,\,\,\,\,\,de\ldots f}^{ab\ldots c}
\end{equation}
where $\left|\frac{\partial x}{\partial\bar{x}}\right|$ is the Jacobian
of the transformation between the two systems. When $w=0$ the tensor
is described as an absolute or true tensor, while when $w=-1$ the
tensor is described as a pseudo tensor. When $w=1$ the tensor may
be described as a tensor density.\footnote{Some of these labels are used differently by different authors.}

$\bullet$ As indicated earlier, a tensor of $m$ contravariant indices
and $n$ covariant indices may be called type ($m,n$). This may be
generalized to include the weight as a third entry and hence the type
of the tensor is identified by ($m,n,w$).

$\bullet$ Relative tensors can be added and subtracted if they are
of the same variance type and have the same weight; the result is
a tensor of the same type and weight. Also, relative tensors can be
equated if they are of the same type and weight.

$\bullet$ Multiplication of relative tensors produces a relative
tensor whose weight is the sum of the weights of the original tensors.
Hence, if the weights are added up to a non-zero value the result
is a relative tensor of that weight; otherwise it is an absolute tensor.

\subsubsection{Isotropic and Anisotropic Tensors}

$\bullet$ Isotropic tensors are characterized by the property that
the values of their components are invariant under coordinate transformation
by proper rotation of axes. In contrast, the values of the components
of anisotropic tensors are dependent on the orientation of the coordinate
axes. Notable examples of isotropic tensors are scalars (rank-0),
the vector $\mathbf{0}$ (rank-1), Kronecker delta $\delta_{ij}$
(rank-2) and Levi-Civita tensor $\epsilon_{ijk}$ (rank-3). Many tensors
describing physical properties of materials, such as stress and magnetic
susceptibility, are anisotropic.

$\bullet$ Direct and inner products of isotropic tensors are isotropic
tensors.

$\bullet$ The zero tensor of any rank is isotropic; therefore if
the components of a tensor vanish in a particular coordinate system
they will vanish in all properly and improperly rotated coordinate
systems.\footnote{For improper rotation, this is more general than being isotropic.}
Consequently, if the components of two tensors are identical in a
particular coordinate system they are identical in all transformed
coordinate systems.

$\bullet$ As indicated, all rank-0 tensors (scalars) are isotropic.
Also, the zero vector, $\mathbf{0}$, of any dimension is isotropic;
in fact it is the only rank-1 isotropic tensor.

\subsubsection{Symmetric and Anti-symmetric Tensors}

$\bullet$ These types of tensor apply to high ranks only (rank $\ge2$).
Moreover, these types are not exhaustive, even for tensors of ranks
$\ge2$, as there are high-rank tensors which are neither symmetric
nor anti-symmetric.

$\bullet$ A rank-2 tensor $A_{ij}$ is symmetric \textit{iff }for
all $i$ and $j$
\begin{equation}
A_{ji}=A_{ij}
\end{equation}
and anti-symmetric or skew-symmetric \textit{iff}
\begin{equation}
A_{ji}=-A_{ij}
\end{equation}
Similar conditions apply to contravariant type tensors (refer also
to the following).

$\bullet$ A rank-$n$ tensor $A_{i_{1}\ldots i_{n}}$ is symmetric
in its two indices $i_{j}$ and $i_{l}$ \textit{iff}
\begin{equation}
A_{i_{1}\ldots i_{l}\ldots i_{j}\ldots i_{n}}=A_{i_{1}\ldots i_{j}\ldots i_{l}\ldots i_{n}}
\end{equation}
and anti-symmetric or skew-symmetric in its two indices $i_{j}$ and
$i_{l}$ \textit{iff}
\begin{equation}
A_{i_{1}\ldots i_{l}\ldots i_{j}\ldots i_{n}}=-A_{i_{1}\ldots i_{j}\ldots i_{l}\ldots i_{n}}
\end{equation}

$\bullet$ Any rank-2 tensor $A_{ij}$ can be synthesized from (or
decomposed into) a symmetric part $A_{(ij)}$ (marked with round brackets
enclosing the indices) and an anti-symmetric part $A_{[ij]}$ (marked
with square brackets) where
\begin{equation}
A_{ij}=A_{(ij)}+A_{[ij]},\,\,\,\,\,\,\,\,\,\,A_{(ij)}=\frac{1}{2}\left(A_{ij}+A_{ji}\right)\,\,\,\,\,\,\,\,\,\&\,\,\,\,\,\,\,\,\,A_{[ij]}=\frac{1}{2}\left(A_{ij}-A_{ji}\right)
\end{equation}

$\bullet$ A rank-3 tensor $A_{ijk}$ can be symmetrized by
\begin{equation}
A_{(ijk)}=\frac{1}{3!}\left(A_{ijk}+A_{kij}+A_{jki}+A_{ikj}+A_{jik}+A_{kji}\right)
\end{equation}
and anti-symmetrized by
\begin{equation}
A_{[ijk]}=\frac{1}{3!}\left(A_{ijk}+A_{kij}+A_{jki}-A_{ikj}-A_{jik}-A_{kji}\right)
\end{equation}

$\bullet$ A rank-$n$ tensor $A_{i_{1}\ldots i_{n}}$ can be symmetrized
by
\begin{equation}
A_{(i_{1}\ldots i_{n})}=\frac{1}{n!}\left(\text{sum of all even \& odd permutations of indices \ensuremath{i}'s}\right)
\end{equation}
and anti-symmetrized by
\begin{equation}
A_{\left[i_{1}\ldots i_{n}\right]}=\frac{1}{n!}\left(\text{sum of all even permutations minus sum of all odd permutations}\right)
\end{equation}

$\bullet$ For a symmetric tensor $\ensuremath{A_{ij}}$ and an anti-symmetric
tensor $\ensuremath{B^{ij}}$ (or the other way around) we have
\begin{equation}
A_{ij}B^{ij}=0
\end{equation}

$\bullet$ The indices whose exchange defines the symmetry and anti-symmetry
relations should be of the same variance type, i.e. both upper or
both lower.

$\bullet$ The symmetry and anti-symmetry characteristic of a tensor
is invariant under coordinate transformation.

$\bullet$ A tensor of high rank ($>2$) may be symmetrized or anti-symmetrized
with respect to only some of its indices instead of all of its indices,
e.g.
\begin{equation}
A_{(ij)k}=\frac{1}{2}\left(A_{ijk}+A_{jik}\right)\,\,\,\,\,\,\,\,\,\,\,\,\,\,\,\&\,\,\,\,\,\,\,\,\,\,\,\,\,\,\,A_{[ij]k}=\frac{1}{2}\left(A_{ijk}-A_{jik}\right)
\end{equation}

$\bullet$ A tensor is totally symmetric \textit{iff}
\begin{equation}
A_{i_{1}\ldots i_{n}}=A_{(i_{1}\ldots i_{n})}
\end{equation}
and totally anti-symmetric \textit{iff}
\begin{equation}
A_{i_{1}\ldots i_{n}}=A_{[i_{1}\ldots i_{n}]}
\end{equation}

$\bullet$ For a totally skew-symmetric tensor (i.e. anti-symmetric
in all of its indices), nonzero entries can occur only when all the
indices are different.

\subsection{Tensor Operations}

$\bullet$ There are many operations that can be performed on tensors
to produce other tensors in general. Some examples of these operations
are addition/subtraction, multiplication by a scalar (rank-0 tensor),
multiplication of tensors (each of rank $>0$), contraction and permutation.
Some of these operations, such as addition and multiplication, involve
more than one tensor while others are performed on a single tensor,
such as contraction and permutation.

$\bullet$ In tensor algebra, division is allowed only for scalars,
hence if the components of an indexed tensor should appear in a denominator,
the tensor should be redefined to avoid this, e.g. $B_{i}=\frac{1}{A_{i}}$.

\subsubsection{Addition and Subtraction}

$\bullet$ Tensors of the same rank and type (covariant/contravariant/mixed
and true/pseudo) can be added algebraically to produce a tensor of
the same rank and type, e.g.
\begin{equation}
a=b+c
\end{equation}
\begin{equation}
A_{i}=B_{i}-C_{i}
\end{equation}
\begin{equation}
A_{j}^{i}=B_{j}^{i}+C_{j}^{i}
\end{equation}

$\bullet$ The added/subtracted terms should have the same indicial
structure with regard to their free indices, as explained in $\S$
\ref{subGeneralRules}, hence $A_{jk}^{i}$ and $B_{ik}^{j}$ cannot
be added or subtracted although they are of the same rank and type,
but $A_{mjk}^{mi}$ and $B_{jk}^{i}$ can be added and subtracted.

$\bullet$ Addition of tensors is associative and commutative:
\begin{equation}
\left(\mathbf{A}+\mathbf{B}\right)+\mathbf{C}=\mathbf{A}+\left(\mathbf{B}+\mathbf{C}\right)
\end{equation}
\begin{equation}
\mathbf{A}+\mathbf{B}=\mathbf{B}+\mathbf{A}
\end{equation}

\subsubsection{Multiplication by Scalar\label{subMultiplicationbyScalar}}

$\bullet$ A tensor can be multiplied by a scalar, which generally
should not be zero, to produce a tensor of the same variance type
and rank, e.g.
\begin{equation}
A_{ik}^{j}=aB_{ik}^{j}
\end{equation}
where $a$ is a non-zero scalar.

$\bullet$ As indicated above, multiplying a tensor by a scalar means
multiplying each component of the tensor by that scalar.

$\bullet$ Multiplication by a scalar is commutative, and associative
when more than two factors are involved.

\subsubsection{Tensor Multiplication \label{subTensorMultiplication} }

$\bullet$ This may also be called outer or exterior or direct or
dyadic multiplication, although some of these names may be reserved
for operations on vectors.

$\bullet$ On multiplying each component of a tensor of rank $r$
by each component of a tensor of rank $k$, both of dimension $m$,
a tensor of rank $(r+k)$ with $m^{r+k}$ components is obtained where
the variance type of each index (covariant or contravariant) is preserved,
e.g.
\begin{equation}
A_{i}B_{j}=C_{ij}
\end{equation}
\begin{equation}
A^{ij}B_{kl}=C_{\,\,\,kl}^{ij}
\end{equation}

$\bullet$ The outer product of a tensor of type ($m,n$) by a tensor
of type ($p,q$) results in a tensor of type ($m+p,n+q$).

$\bullet$ Direct multiplication of tensors may be marked by the symbol
$\varotimes$, mostly when using symbolic notation for tensors, e.g.
$\mathbf{A}\varotimes\mathbf{B}$. However, in the present text no
symbol will be used for the operation of direct multiplication.

$\bullet$ Direct multiplication of tensors is not commutative.

$\bullet$ The outer product operation is distributive with respect
to the algebraic sum of tensors:
\begin{equation}
\mathbf{A}\left(\mathbf{B}\pm\mathbf{C}\right)=\mathbf{A}\mathbf{B}\pm\mathbf{A}\mathbf{C}\,\,\,\,\,\,\,\,\,\,\,\&\,\,\,\,\,\,\,\,\,\,\,\,\left(\mathbf{B}\pm\mathbf{C}\right)\mathbf{A}=\mathbf{B}\mathbf{A}\pm\mathbf{C\mathbf{A}}
\end{equation}

$\bullet$ Multiplication of a tensor by a scalar (refer to $\S$
\ref{subMultiplicationbyScalar}) may be regarded as a special case
of direct multiplication.

$\bullet$ The rank-2 tensor constructed as a result of the direct
multiplication of two vectors is commonly called dyad.

$\bullet$ Tensors may be expressed as an outer product of vectors
where the rank of the resultant product is equal to the number of
the vectors involved (e.g. 2 for dyads and 3 for triads).

$\bullet$ Not every tensor can be synthesized as a product of lower
rank tensors.

$\bullet$ In the outer product, it is understood that all the indices
of the involved tensors have the same range.

$\bullet$ The outer product of tensors yields a tensor.

\subsubsection{Contraction\label{subsubContraction}}

$\bullet$ Contraction of a tensor of rank $>1$ is to make two free
indices identical, by unifying their symbols, and perform summation
over these repeated indices, e.g.
\begin{equation}
A_{i}^{j}\,\,\,\,\,\,\,\,\,\,\underrightarrow{\mathrm{contraction\,\,}}\,\,\,\,\,\,\,\,\,\,A_{i}^{i}
\end{equation}
\begin{equation}
A_{il}^{jk}\,\,\,\,\,\,\,\,\,\,\underrightarrow{\mathrm{contraction\,\,on}\,\,jl\,\,}\,\,\,\,\,\,\,\,\,\,A_{im}^{mk}
\end{equation}

$\bullet$ Contraction results in a reduction of the rank by 2 since
it implies the annihilation of two free indices. Therefore, the contraction
of a rank-2 tensor is a scalar, the contraction of a rank-3 tensor
is a vector, the contraction of a rank-4 tensor is a rank-2 tensor,
and so on.

$\bullet$ For general non-Cartesian coordinate systems, the pair
of contracted indices should be different in their variance type,
i.e. one upper and one lower. Hence, contraction of a mixed tensor
of type ($m,n$) will, in general, produce a tensor of type ($m-1,n-1$).

$\bullet$ A tensor of type ($p,q$) can have $p\times q$ possible
contractions, i.e. one contraction for each pair of lower and upper
indices.

$\bullet$ A common example of contraction is the dot product operation
on vectors which can be regarded as a direct multiplication (refer
to $\S$ \ref{subTensorMultiplication}) of the two vectors, which
results in a rank-2 tensor, followed by a contraction.

$\bullet$ In matrix algebra, taking the trace (summing the diagonal
elements) can also be considered as contraction of the matrix, which
under certain conditions can represent a rank-2 tensor, and hence
it yields the trace which is a scalar.

$\bullet$ Applying the index contraction operation on a tensor results
into a tensor.

$\bullet$ Application of contraction of indices operation on a relative
tensor (see $\S$ \ref{subAbsoluteRelative}) produces a relative
tensor of the same weight as the original tensor.

\subsubsection{Inner Product\label{secInnerProduct}}

$\bullet$ On taking the outer product (refer to $\S$ \ref{subTensorMultiplication})
of two tensors of rank $\ge1$ followed by a contraction on two indices
of the product, an inner product of the two tensors is formed. Hence
if one of the original tensors is of rank-$m$ and the other is of
rank-$n$, the inner product will be of rank-($m+n-2$).

$\bullet$ The inner product operation is usually symbolized by a
single dot between the two tensors, e.g. $\mathbf{A}\cdot\mathbf{B}$,
to indicate contraction following outer multiplication.

$\bullet$ In general, the inner product is not commutative. When
one or both of the tensors involved in the inner product are of rank
$>1$ the order of the multiplicands does matter.

$\bullet$ The inner product operation is distributive with respect
to the algebraic sum of tensors:
\begin{equation}
\mathbf{A}\cdot\left(\mathbf{B}\pm\mathbf{C}\right)=\mathbf{A}\cdot\mathbf{B}\pm\mathbf{A}\cdot\mathbf{C}\,\,\,\,\,\,\,\,\,\,\,\&\,\,\,\,\,\,\,\,\,\,\,\,\left(\mathbf{B}\pm\mathbf{C}\right)\cdot\mathbf{A}=\mathbf{B}\cdot\mathbf{A}\pm\mathbf{C}\cdot\mathbf{A}
\end{equation}

$\bullet$ As indicated before (see $\S$ \ref{subsubContraction}),
the dot product of two vectors is an example of the inner product
of tensors, i.e. it is an inner product of two rank-1 tensors to produce
a rank-0 tensor:
\begin{equation}
\left[\mathbf{a}\mathbf{b}\right]_{i}^{\,\,j}=a_{i}b^{j}\,\,\,\,\,\,\,\,\,\,\underrightarrow{\mathrm{contraction\,\,}}\,\,\,\,\,\,\,\,\,\,\mathbf{a}\cdot\mathbf{b}=a_{i}b^{i}
\end{equation}

$\bullet$ Another common example (from linear algebra) of inner product
is the multiplication of a matrix (representing a rank-2 tensor assuming
certain conditions) by a vector (rank-1 tensor) to produce a vector,
e.g.
\begin{equation}
\left[\mathbf{A}\mathbf{b}\right]_{ij}^{\,\,\,k}=A_{ij}b^{k}\,\,\,\,\,\,\,\,\,\,\underrightarrow{\mathrm{contraction\,\,on}\,\,jk\,\,}\,\,\,\,\,\,\,\,\,\,\left[\mathbf{A}\cdot\mathbf{b}\right]_{i}=A_{ij}b^{j}
\end{equation}

The multiplication of two $n\times n$ matrices is another example
of inner product (see Eq. \ref{eqMatrixMultiplication}).

$\bullet$ For tensors whose outer product produces a tensor of rank
$>2$, various contraction operations between different sets of indices
can occur and hence more than one inner product, which are different
in general, can be defined. Moreover, when the outer product produces
a tensor of rank $>3$ more than one contraction can take place simultaneously.

$\bullet$ There are more specialized types of inner product; some
of which may be defined differently by different authors. For example,
a double inner product of two rank-2 tensors, $\mathbf{A}$ and $\mathbf{B}$,
may be defined and denoted by double vertically- or horizontally-aligned
dots (e.g. $\mathbf{A}\colon\mathbf{B}$ or $\mathbf{A}\cdot\cdot\mathbf{B}$)
to indicate double contraction taking place between different pairs
of indices. An instance of these types is the inner product with double
contraction of two dyads which is commonly defined by\footnote{It is also defined differently by some authors.}
\begin{equation}
\mathbf{a}\mathbf{b}\colon\mathbf{c}\mathbf{d}=\left(\mathbf{a}\cdot\mathbf{c}\right)\left(\mathbf{b}\cdot\mathbf{d}\right)
\end{equation}
and the result is a scalar. The single dots in the right hand side
of the last equation symbolize the conventional dot product of two
vectors. Some authors may define a different type of double-contraction
inner product of two dyads, symbolized by two horizontally-aligned
dots, which may be called a ``transposed contraction'', and is given
by
\begin{equation}
\mathbf{a}\mathbf{b}\cdot\cdot\mathbf{c}\mathbf{d}=\mathbf{a}\mathbf{b}\colon\mathbf{d}\mathbf{c}=\left(\mathbf{a}\cdot\mathbf{d}\right)\left(\mathbf{b}\cdot\mathbf{c}\right)
\end{equation}
where the result is also a scalar. However, different authors may
have different conventions and hence one should be vigilant about
such differences.

$\bullet$ For two rank-2 tensors, the aforementioned double-contraction
inner products are similarly defined as in the case of two dyads:
\begin{equation}
\mathbf{A}\colon\mathbf{B}=A_{ij}B_{ij}\,\,\,\,\,\,\,\,\,\,\&\,\,\,\,\,\,\,\,\,\,\mathbf{A}\cdot\cdot\mathbf{B}=A_{ij}B_{ji}
\end{equation}

$\bullet$ Inner products with higher multiplicity of contractions
are similarly defined, and hence can be regarded as trivial extensions
of the inner products with lower contraction multiplicities.

$\bullet$ The inner product of tensors produces a tensor because
the inner product is an outer product operation followed by an index
contraction operation and both of these operations on tensors produce
tensors.

\subsubsection{Permutation}

$\bullet$ A tensor may be obtained by exchanging the indices of another
tensor, e.g. transposition of rank-2 tensors.

$\bullet$ Tensor permutation applies only to tensors of rank $\ge2$.

$\bullet$ The collection of tensors obtained by permuting the indices
of a basic tensor may be called isomers.

\subsection{Tensor Test: Quotient Rule}

$\bullet$ Sometimes a tensor-like object may be suspected for being
a tensor; in such cases a test based on the ``quotient rule'' can
be used to clarify the situation. According to this rule, if the inner
product of a suspected tensor with a known tensor is a tensor then
the suspect is a tensor. In more formal terms, if it is not known
if $\mathbf{A}$ is a tensor but it is known that $\mathbf{B}$ and
$\mathbf{C}$ are tensors; moreover it is known that the following
relation holds true in all rotated (properly-transformed) coordinate
frames:
\begin{equation}
A_{pq\ldots k\ldots m}B_{ij\ldots k\ldots n}=C_{pq\ldots mij\ldots n}\label{eqQuotientRule}
\end{equation}
then $\mathbf{A}$ is a tensor. Here, $\mathbf{A}$, $\mathbf{B}$
and $\mathbf{C}$ are respectively of ranks $m,\,n$ and ($m+n-2$),
due to the contraction on $k$ which can be any index of $\mathbf{A}$
and $\mathbf{B}$ independently.

$\bullet$ Testing for being a tensor can also be done by applying
first principles through direct substitution in the transformation
equations. However, using the quotient rule is generally more convenient
and requires less work.

$\bullet$ The quotient rule may be considered as a replacement for
the division operation which is not defined for tensors.

\pagebreak{}

\section{$\mathbf{\delta}$ and $\mathbf{\epsilon}$ Tensors}

$\bullet$ These tensors are of particular importance in tensor calculus
due to their distinctive properties and unique transformation attributes.
They are numerical tensors with fixed components in all coordinate
systems. The first is called Kronecker delta or unit tensor, while
the second is called Levi-Civita\footnote{This name is usually used for the rank-3 tensor. Also some authors
distinguish between the permutation tensor and Levi-Civita tensor
even for rank-3. Moreover, some of the common labels and descriptions
of $\epsilon$ are more specific to rank-3.}, permutation, anti-symmetric and alternating tensor.

$\bullet$ The $\delta$ and $\epsilon$ tensors are conserved under
coordinate transformations and hence they are the same for all systems
of coordinate.\footnote{For the permutation tensor, the statement applies to proper coordinate
transformations.}

\subsection{Kronecker $\delta$\label{subKronecker}}

$\bullet$ This is a rank-2 symmetric tensor in all dimensions, i.e.
\begin{equation}
\delta_{ij}=\delta_{ji}\,\,\,\,\,\,\,\,\,\,\,\,\,\left(i,j=1,2,\ldots,n\right)
\end{equation}
Similar identities apply to the contravariant and mixed types of this
tensor.

$\bullet$ It is invariant in all coordinate systems, and hence it
is an isotropic tensor.\footnote{In fact it is more general than isotropic as it is invariant even
under improper coordinate transformations.}

$\bullet$ It is defined as:
\begin{equation}
\delta_{ij}=\begin{cases}
1 & (i=j)\\
0\,\,\,\,\,\,\,\,\,\,\,\,\,\, & (i\neq j)
\end{cases}\label{eqKroneckerDefinitionNormal}
\end{equation}
and hence it can be considered as the identity matrix, e.g. for 3D
\begin{equation}
\left[\delta_{ij}\right]=\left[\begin{array}{ccc}
\delta_{11} & \delta_{12} & \delta_{13}\\
\delta_{21} & \delta_{22} & \delta_{23}\\
\delta_{31} & \delta_{32} & \delta_{33}
\end{array}\right]=\left[\begin{array}{ccc}
1 & 0 & 0\\
0 & 1 & 0\\
0 & 0 & 1
\end{array}\right]
\end{equation}

$\bullet$ Covariant, contravariant and mixed type of this tensor
are the same, that is
\begin{equation}
\delta_{\,j}^{i}=\delta_{i}^{\,j}=\delta^{ij}=\delta_{ij}
\end{equation}

\subsection{Permutation $\epsilon$\label{subPermutation}}

$\bullet$ This is an isotropic tensor. It has a rank equal to the
number of dimensions; hence, a rank-$n$ permutation tensor has $n^{n}$
components.

$\bullet$ It is totally anti-symmetric in each pair of its indices,
i.e. it changes sign on swapping any two of its indices, that is
\begin{equation}
\epsilon_{i_{1}\ldots i_{k}\ldots i_{l}\ldots i_{n}}=-\epsilon_{i_{1}\ldots i_{l}\ldots i_{k}\ldots i_{n}}
\end{equation}
The reason is that any exchange of two indices requires an even/odd
number of single-step shifts to the right of the first index plus
an odd/even number of single-step shifts to the left of the second
index, so the total number of shifts is odd and hence it is an odd
permutation of the original arrangement.

$\bullet$ It is a pseudo tensor since it acquires a minus sign under
improper orthogonal transformation of coordinates (inversion of axes
with possible superposition of rotation).

$\bullet$ Definition of rank-2 $\epsilon$ ($\epsilon_{ij}$):
\begin{equation}
\epsilon_{12}=1,\,\,\,\,\,\,\,\,\,\,\epsilon_{21}=-1\,\,\,\,\,\,\,\,\,\,\&\,\,\,\,\,\,\,\,\,\,\epsilon_{11}=\epsilon_{22}=0
\end{equation}

$\bullet$ Definition of rank-3 $\epsilon$ ($\epsilon_{ijk}$):
\begin{equation}
\epsilon_{ijk}=\begin{cases}
\,\,\,\,\,1 & (i,j,k\text{ is even permutation of 1,2,3})\\
-1 & (i,j,k\text{ is odd permutation of 1,2,3})\\
\,\,\,\,\,0\,\,\,\,\,\,\,\,\,\,\,\,\,\, & (\text{repeated\,index})
\end{cases}\label{eqEpsilon3Definition}
\end{equation}

$\bullet$ The definition of rank-$n$ $\epsilon$ ($\epsilon_{i_{1}i_{2}\ldots i_{n}}$)
is similar to the definition of rank-3 $\epsilon$ considering index
repetition and even or odd permutations of its indices $\left(i_{1},i_{2},\cdots,i_{n}\right)$
corresponding to $\left(1,2,\cdots,n\right)$, that is
\begin{equation}
\epsilon_{i_{1}i_{2}\ldots i_{n}}=\begin{cases}
\,\,\,\,\,1 & \left[\left(i_{1},i_{2},\ldots,i_{n}\right)\text{ is even permutation of (\ensuremath{1,2,\ldots,n})}\right]\\
-1 & \left[\left(i_{1},i_{2},\ldots,i_{n}\right)\text{ is odd permutation of (\ensuremath{1,2,\ldots,n})}\right]\\
\,\,\,\,\,0\,\,\,\,\,\,\,\,\,\,\,\,\,\, & \text{\ensuremath{\left[\mathrm{repeated\,index}\right]}}
\end{cases}\label{eqEpsilonnDefinition}
\end{equation}

$\bullet$ $\epsilon$ may be considered a contravariant relative
tensor of weight $+1$ or a covariant relative tensor of weight $-1$.
Hence, in 2D, 3D and $n$D spaces respectively we have:
\begin{eqnarray}
\epsilon_{ij} & = & \epsilon^{ij}\\
\epsilon_{ijk} & = & \epsilon^{ijk}\\
\epsilon_{i_{1}i_{2}\ldots i_{n}} & = & \epsilon^{i_{1}i_{2}\ldots i_{n}}
\end{eqnarray}

\subsection{Useful Identities Involving $\delta$ or/and $\epsilon$}

\subsubsection{Identities Involving $\delta$}

$\bullet$ When an index of the Kronecker delta is involved in a contraction
operation by repeating an index in another tensor in its own term,
the effect of this is to replace the shared index in the other tensor
by the other index of the Kronecker delta, that is
\begin{equation}
\delta_{ij}A_{j}=A_{i}\label{EqIndexReplace}
\end{equation}
In such cases the Kronecker delta is described as the substitution
or index replacement operator. Hence,
\begin{equation}
\delta_{ij}\delta_{jk}=\delta_{ik}
\end{equation}
Similarly,
\begin{equation}
\delta_{ij}\delta_{jk}\delta_{ki}=\delta_{ik}\delta_{ki}=\delta_{ii}=n\label{eqdeltas}
\end{equation}
where $n$ is the space dimension.

$\bullet$ Because the coordinates are independent of each other:
\begin{equation}
\frac{\partial x_{i}}{\partial x_{j}}=\partial_{j}x_{i}=x_{i,j}=\delta_{ij}\label{eqdxdelta}
\end{equation}
Hence, in an $n$D space we have
\begin{equation}
\partial_{i}x_{i}=\delta_{ii}=n\label{eqdxn}
\end{equation}

$\bullet$ For orthonormal Cartesian systems:
\begin{equation}
\frac{\partial x^{i}}{\partial x^{j}}=\frac{\partial x^{j}}{\partial x^{i}}=\delta_{ij}=\delta^{ij}\label{eqdxdxdelta}
\end{equation}

$\bullet$ For a set of orthonormal basis vectors in orthonormal Cartesian
systems:
\begin{equation}
\mathbf{e}_{i}\cdot\mathbf{e}_{j}=\delta_{ij}
\end{equation}

$\bullet$ The double inner product of two dyads formed by orthonormal
basis vectors of an orthonormal Cartesian system is given by:
\begin{equation}
\mathbf{e}_{i}\mathbf{e}_{j}\colon\mathbf{e}_{k}\mathbf{e}_{l}=\delta_{ik}\delta_{jl}
\end{equation}

\subsubsection{Identities Involving $\epsilon$}

$\bullet$ For rank-3 $\epsilon$:
\begin{equation}
\epsilon_{ijk}=\epsilon_{kij}=\epsilon_{jki}=-\epsilon_{ikj}=-\epsilon_{jik}=-\epsilon_{kji}\,\,\,\,\,\,\,\,\,\,\text{(sense of cyclic order)}\label{EqEpsilonCycle}
\end{equation}
These equations demonstrate the fact that rank-3 $\epsilon$ is totally
anti-symmetric in all of its indices since a shift of any two indices
reverses the sign. This also reflects the fact that the above tensor
system has only one independent component.

$\bullet$ For rank-2 $\epsilon$:
\begin{equation}
\epsilon_{ij}=\left(j-i\right)
\end{equation}

$\bullet$ For rank-3 $\epsilon$:
\begin{equation}
\epsilon_{ijk}=\frac{1}{2}\left(j-i\right)\left(k-i\right)\left(k-j\right)
\end{equation}

$\bullet$ For rank-4 $\epsilon$:
\begin{equation}
\epsilon_{ijkl}=\frac{1}{12}\left(j-i\right)\left(k-i\right)\left(l-i\right)\left(k-j\right)\left(l-j\right)\left(l-k\right)
\end{equation}

$\bullet$ For rank-$n$ $\epsilon$:
\begin{equation}
\epsilon_{a_{1}a_{2}\cdots a_{n}}=\prod_{i=1}^{n-1}\left[\frac{1}{i!}\prod_{j=i+1}^{n}\left(a_{j}-a_{i}\right)\right]=\frac{1}{S(n-1)}\prod_{1\le i<j\le n}\left(a_{j}-a_{i}\right)
\end{equation}
where $S(n-1)$ is the super-factorial function of $(n-1)$ which
is defined as
\begin{equation}
S(k)=\prod_{i=1}^{k}i!=1!\cdot2!\cdot\ldots\cdot k!
\end{equation}
A simpler formula for rank-$n$ $\epsilon$ can be obtained from the
previous one by ignoring the magnitude of the multiplication factors
and taking only their signs, that is
\begin{equation}
\epsilon_{a_{1}a_{2}\cdots a_{n}}=\prod_{1\le i<j\le n}\sigma\left(a_{j}-a_{i}\right)=\sigma\left(\prod_{1\le i<j\le n}\left(a_{j}-a_{i}\right)\right)
\end{equation}
where
\begin{equation}
\sigma(k)=\begin{cases}
+1 & (k>0)\\
-1 & (k<0)\\
\,\,\,\,\,0\,\,\,\,\,\,\,\,\,\,\,\,\,\, & (k=0)
\end{cases}
\end{equation}

$\bullet$ For rank-$n$ $\epsilon$:
\begin{equation}
\epsilon_{i_{1}i_{2}\cdots i_{n}}\,\epsilon_{i_{1}i_{2}\cdots i_{n}}=n!
\end{equation}
because this is the sum of the squares of $\epsilon_{i_{1}i_{2}\cdots i_{n}}$
over all the permutations of $n$ different indices which is equal
to $n!$ where the value of $\epsilon$ of each one of these permutations
is either $+1$ or $-1$ and hence in both cases their square is 1.

$\bullet$ For a symmetric tensor $A_{jk}$:
\begin{equation}
\epsilon_{ijk}A_{jk}=0
\end{equation}
because an exchange of the two indices of $A_{jk}$ does not affect
its value due to the symmetry whereas a similar exchange in these
indices in $\epsilon_{ijk}$ results in a sign change; hence each
term in the sum has its own negative and therefore the total sum will
vanish.

$\bullet$
\begin{equation}
\epsilon_{ijk}A_{i}A_{j}=\epsilon_{ijk}A_{i}A_{k}=\epsilon_{ijk}A_{j}A_{k}=0\label{eqPermutingTwoFactors}
\end{equation}
because, due to the commutativity of multiplication, an exchange of
the indices in $A$'s will not affect the value but a similar exchange
in the corresponding indices of $\epsilon_{ijk}$ will cause a change
in sign; hence each term in the sum has its own negative and therefore
the total sum will be zero.

$\bullet$ For a set of orthonormal basis vectors in a 3D space with
a right-handed orthonormal Cartesian coordinate system:
\begin{equation}
\mathbf{e}_{i}\times\mathbf{e}_{j}=\epsilon_{ijk}\mathbf{e}_{k}
\end{equation}
\begin{equation}
\mathbf{e}_{i}\cdot\left(\mathbf{e}_{j}\times\mathbf{e}_{k}\right)=\epsilon_{ijk}
\end{equation}

\subsubsection{Identities Involving $\delta$ and $\epsilon$}

$\bullet$
\begin{equation}
\epsilon_{ijk}\delta_{1i}\delta_{2j}\delta_{3k}=\epsilon_{123}=1
\end{equation}

$\bullet$ For rank-2 $\epsilon$:
\begin{equation}
\epsilon_{ij}\epsilon_{kl}=\begin{vmatrix}\begin{array}{cc}
\delta_{ik} & \delta_{il}\\
\delta_{jk} & \delta_{jl}
\end{array}\end{vmatrix}=\delta_{ik}\delta_{jl}-\delta_{il}\delta_{jk}
\end{equation}
\begin{equation}
\epsilon_{il}\epsilon_{kl}=\delta_{ik}
\end{equation}
\begin{equation}
\epsilon_{ij}\epsilon_{ij}=2
\end{equation}

$\bullet$ For rank-3 $\epsilon$:
\begin{equation}
\epsilon_{ijk}\epsilon_{lmn}=\begin{vmatrix}\begin{array}{ccc}
\delta_{il} & \delta_{im} & \delta_{in}\\
\delta_{jl} & \delta_{jm} & \delta_{jn}\\
\delta_{kl} & \delta_{km} & \delta_{kn}
\end{array}\end{vmatrix}=\delta_{il}\delta_{jm}\delta_{kn}+\delta_{im}\delta_{jn}\delta_{kl}+\delta_{in}\delta_{jl}\delta_{km}-\delta_{il}\delta_{jn}\delta_{km}-\delta_{im}\delta_{jl}\delta_{kn}-\delta_{in}\delta_{jm}\delta_{kl}
\end{equation}
\begin{equation}
\epsilon_{ijk}\epsilon_{lmk}=\begin{vmatrix}\begin{array}{cc}
\delta_{il} & \delta_{im}\\
\delta_{jl} & \delta_{jm}
\end{array}\end{vmatrix}=\delta_{il}\delta_{jm}-\delta_{im}\delta_{jl}\label{EqEpsilonDelta}
\end{equation}
The last identity is very useful in manipulating and simplifying tensor
expressions and proving vector and tensor identities.
\begin{equation}
\epsilon_{ijk}\epsilon_{ljk}=2\delta_{il}
\end{equation}
\begin{equation}
\epsilon_{ijk}\epsilon_{ijk}=2\delta_{ii}=6
\end{equation}
since the rank and dimension of $\epsilon$ are the same, which is
3 in this case.

$\bullet$ For rank-$n$ $\epsilon$:
\begin{equation}
\epsilon_{i_{1}i_{2}\cdots i_{n}}\,\epsilon_{j_{1}j_{2}\cdots j_{n}}=\begin{vmatrix}\begin{array}{cccc}
\delta_{i_{1}j_{1}} & \delta_{i_{1}j_{2}} & \cdots & \delta_{i_{1}j_{n}}\\
\delta_{i_{2}j_{1}} & \delta_{i_{2}j_{2}} & \cdots & \delta_{i_{2}j_{n}}\\
\vdots & \vdots & \ddots & \vdots\\
\delta_{i_{n}j_{1}} & \delta_{i_{n}j_{2}} & \cdots & \delta_{i_{n}j_{n}}
\end{array}\end{vmatrix}\label{eqEpsilon2}
\end{equation}

$\bullet$ According to Eqs. \ref{eqEpsilon3Definition} and \ref{EqIndexReplace}:
\begin{equation}
\epsilon_{ijk}\delta_{ij}=\epsilon_{ijk}\delta_{ik}=\epsilon_{ijk}\delta_{jk}=0
\end{equation}

\subsection{Generalized Kronecker delta}

$\bullet$ The generalized Kronecker delta is defined by:
\begin{equation}
\delta_{j_{1}\ldots j_{n}}^{i_{1}\ldots i_{n}}=\begin{cases}
\,\,\,\,\,1 & \left[(j_{1}\ldots j_{n})\text{ is even permutation of (\ensuremath{i_{1}\ldots i_{n})}}\right]\\
-1 & \left[(j_{1}\ldots j_{n})\text{ is odd permutation of (\ensuremath{i_{1}\ldots i_{n})}}\right]\\
\,\,\,\,\,0\,\,\,\,\,\,\,\,\,\,\,\,\,\, & \left[\text{repeated\,\,\ensuremath{j}'s}\right]
\end{cases}\label{eqGeneralizedKronecker}
\end{equation}
It can also be defined by the following $n\times n$ determinant:
\begin{equation}
\delta_{j_{1}\ldots j_{n}}^{i_{1}\ldots i_{n}}=\begin{vmatrix}\begin{array}{cccc}
\delta_{j_{1}}^{i_{1}} & \delta_{j_{2}}^{i_{1}} & \cdots & \delta_{j_{n}}^{i_{1}}\\
\delta_{j_{1}}^{i_{2}} & \delta_{j_{2}}^{i_{2}} & \cdots & \delta_{j_{n}}^{i_{2}}\\
\vdots & \vdots & \ddots & \vdots\\
\delta_{j_{1}}^{i_{n}} & \delta_{j_{2}}^{i_{n}} & \cdots & \delta_{j_{n}}^{i_{n}}
\end{array}\end{vmatrix}\label{eqGeneralizedKronecker2}
\end{equation}
where the $\delta_{j}^{i}$ entries in the determinant are the normal
Kronecker delta as defined by Eq. \ref{eqKroneckerDefinitionNormal}.

$\bullet$ Accordingly, the relation between the rank-$n$ $\epsilon$
and the generalized Kronecker delta in an $n$D space is given by:
\begin{equation}
\epsilon_{i_{1}i_{2}\ldots i_{n}}=\delta_{i_{1}i_{2}\ldots i_{n}}^{1\,2\ldots n}\,\,\,\,\,\,\,\,\,\,\&\,\,\,\,\,\,\,\,\,\,\epsilon^{i_{1}i_{2}\ldots i_{n}}=\delta_{1\,2\ldots n}^{i_{1}i_{2}\ldots i_{n}}
\end{equation}
Hence, the permutation tensor $\epsilon$ may be considered as a special
case of the generalized Kronecker delta. Consequently the permutation
symbol can be written as an $n\times n$ determinant consisting of
the normal Kronecker deltas.

$\bullet$ If we define
\begin{equation}
\delta_{lm}^{ij}=\delta_{lmk}^{ijk}
\end{equation}
then Eq. \ref{EqEpsilonDelta} will take the following form:
\begin{equation}
\delta_{lm}^{ij}=\delta_{l}^{i}\delta_{m}^{j}-\delta_{m}^{i}\delta_{l}^{j}
\end{equation}
Other identities involving $\delta$ and $\epsilon$ can also be formulated
in terms of the generalized Kronecker delta.

$\bullet$ On comparing Eq. \ref{eqEpsilon2} with Eq. \ref{eqGeneralizedKronecker2}
we conclude
\begin{equation}
\delta_{j_{1}\ldots j_{n}}^{i_{1}\ldots i_{n}}=\epsilon^{i_{1}\ldots i_{n}}\,\epsilon_{j_{1}\ldots j_{n}}
\end{equation}

\pagebreak{}

\section{Applications of Tensor Notation and Techniques}

\subsection{Common Definitions in Tensor Notation\label{subCommonDefinitions}}

$\bullet$ The trace of a matrix $\mathbf{A}$ representing a rank-2
tensor is:
\begin{equation}
\mathrm{tr}\left(\mathbf{A}\right)=A_{ii}
\end{equation}

$\bullet$ For a $3\times3$ matrix representing a rank-2 tensor in
a 3D space, the determinant is:
\begin{equation}
\mathrm{det}\left(\mathbf{A}\right)=\begin{vmatrix}\begin{array}{ccc}
A_{11} & A_{12} & A_{13}\\
A_{21} & A_{22} & A_{23}\\
A_{31} & A_{32} & A_{33}
\end{array}\end{vmatrix}=\epsilon_{ijk}A_{1i}A_{2j}A_{3k}=\epsilon_{ijk}A_{i1}A_{j2}A_{k3}
\end{equation}
where the last two equalities represent the expansion of the determinant
by row and by column. Alternatively
\begin{equation}
\mathrm{det}\left(\mathbf{A}\right)=\frac{1}{3!}\epsilon_{ijk}\epsilon_{lmn}A_{il}A_{jm}A_{kn}
\end{equation}

$\bullet$ For an $n\times n$ matrix representing a rank-2 tensor
in an $n$D space, the determinant is:
\begin{equation}
\mathrm{det}\left(\mathbf{A}\right)=\epsilon_{i_{1}\cdots i_{n}}A_{1i_{1}}\ldots A_{ni_{n}}=\epsilon_{i_{1}\cdots i_{n}}A_{i_{1}1}\ldots A_{i_{n}n}=\frac{1}{n!}\epsilon_{i_{1}\cdots i_{n}}\,\epsilon_{j_{1}\cdots j_{n}}A_{i_{1}j_{1}}\ldots A_{i_{n}j_{n}}
\end{equation}

$\bullet$ The inverse of a matrix $\mathbf{A}$ representing a rank-2
tensor is:
\begin{equation}
\left[\mathbf{A}^{-1}\right]_{ij}=\frac{1}{2\,\mathrm{det}\left(\mathbf{A}\right)}\epsilon_{jmn}\,\epsilon_{ipq}A_{mp}A_{nq}
\end{equation}

$\bullet$ The multiplication of a matrix $\mathbf{A}$ by a vector
$\mathbf{b}$ as defined in linear algebra is:
\begin{equation}
\left[\mathbf{A}\mathbf{b}\right]_{i}=A_{ij}b_{j}
\end{equation}
It should be noticed that here we are using matrix notation. The multiplication
operation, according to the symbolic notation of tensors, should be
denoted by a dot between the tensor and the vector, i.e. $\mathbf{A}\mathbf{\cdot b}$.\footnote{The matrix multiplication in matrix notation is equivalent to a dot
product operation in tensor notation.}

$\bullet$ The multiplication of two $n\times n$ matrices $\mathbf{A}$
and $\mathbf{B}$ as defined in linear algebra is:
\begin{equation}
\left[\mathbf{A}\mathbf{B}\right]_{ik}=A_{ij}B_{jk}\label{eqMatrixMultiplication}
\end{equation}
Again, here we are using matrix notation; otherwise a dot should be
inserted between the two matrices.

$\bullet$ The dot product of two vectors is:
\begin{equation}
\mathbf{A}\cdot\mathbf{B=}\delta_{ij}A_{i}B_{j}=A_{i}B_{i}\label{eqDotProduct}
\end{equation}
The readers are referred to $\S$ \ref{secInnerProduct} for a more
general definition of this type of product that includes higher rank
tensors.

$\bullet$ The cross product of two vectors is:
\begin{equation}
\left[\mathbf{A}\times\mathbf{B}\right]_{i}=\epsilon_{ijk}A_{j}B_{k}\label{EqCrossProduct}
\end{equation}

$\bullet$ The scalar triple product of three vectors is:
\begin{equation}
\mathbf{A}\cdot\left(\mathbf{B}\times\mathbf{C}\right)=\begin{vmatrix}\begin{array}{ccc}
A_{1} & A_{2} & A_{3}\\
B_{1} & B_{2} & B_{3}\\
C_{1} & C_{2} & C_{3}
\end{array}\end{vmatrix}=\epsilon_{ijk}A_{i}B_{j}C_{k}\label{EqScalarTripleProduct}
\end{equation}

$\bullet$ The vector triple product of three vectors is:
\begin{equation}
\left[\mathbf{A}\times\left(\mathbf{B}\times\mathbf{C}\right)\right]_{i}=\epsilon_{ijk}\epsilon_{klm}A_{j}B_{l}C_{m}
\end{equation}

\subsection{Scalar Invariants of Tensors}

$\bullet$ In the following we list and write in tensor notation a
number of invariants of low rank tensors which have special importance
due to their widespread applications in vector and tensor calculus.
All These invariants are scalars.

$\bullet$ The value of a scalar (rank-0 tensor), which consists of
a magnitude and a sign, is invariant under coordinate transformation.

$\bullet$ An invariant of a vector (rank-1 tensor) under coordinate
transformations is its magnitude, i.e. length (the direction is also
invariant but it is not scalar!).\footnote{In fact the magnitude alone is invariant under coordinate transformations
even for pseudo vectors because it is a scalar.}

$\bullet$ The main three independent scalar invariants of a rank-2
tensor $\mathbf{A}$ under change of basis are:
\begin{equation}
I=\mathrm{tr}\left(\mathbf{A}\right)=A_{ii}
\end{equation}
\begin{equation}
II=\mathrm{tr}\left(\mathbf{A}^{2}\right)=A_{ij}A_{ji}
\end{equation}
\begin{equation}
III=\mathrm{tr}\left(\mathbf{A}^{3}\right)=A_{ij}A_{jk}A_{ki}
\end{equation}

$\bullet$ Different forms of the three invariants of a rank-2 tensor
$\mathbf{A}$, which are also widely used, are:
\begin{equation}
I_{1}=I=A_{ii}
\end{equation}
\begin{equation}
I_{2}=\frac{1}{2}\left(I^{2}-II\right)=\frac{1}{2}\left(A_{ii}A_{jj}-A_{ij}A_{ji}\right)
\end{equation}
\begin{equation}
I_{3}=\mathrm{det}\left(\mathbf{A}\right)=\frac{1}{3!}\left(I^{3}-3I\,\,II+2III\right)=\frac{1}{3!}\epsilon_{ijk}\epsilon_{pqr}A_{ip}A_{jq}A_{kr}
\end{equation}

$\bullet$ The invariants $I$, $II$ and $III$ can similarly be
defined in terms of the invariants $I_{1}$, $I_{2}$ and $I_{3}$
as follow:
\begin{equation}
I=I_{1}
\end{equation}
\begin{equation}
II=I_{1}^{2}-2I_{2}
\end{equation}
\begin{equation}
III=I_{1}^{3}-3I_{1}I_{2}+3I_{3}
\end{equation}

$\bullet$ Since the determinant of a matrix representing a rank-2
tensor is invariant, then if the determinant vanishes in one coordinate
system it will vanish in all coordinate systems and vice versa. Consequently,
if a rank-2 tensor is invertible in a particular coordinate system,
it is invertible in all coordinate systems.

$\bullet$ Ten joint invariants between two rank-2 tensors, $\mathbf{A}$
and $\mathbf{B}$, can be formed; these are: $\mathrm{tr}\left(\mathbf{A}\right)$,
$\mathrm{tr}\left(\mathbf{B}\right)$, $\mathrm{tr}\left(\mathbf{A}^{2}\right)$,
$\mathrm{tr}\left(\mathbf{B}^{2}\right)$, $\mathrm{tr}\left(\mathbf{A}^{3}\right)$,
$\mathrm{tr}\left(\mathbf{B}^{3}\right)$, $\mathrm{tr}\left(\mathbf{A}\cdot\mathbf{B}\right)$,
$\mathrm{tr}\left(\mathbf{A}^{2}\cdot\mathbf{B}\right)$, $\mathrm{tr}\left(\mathbf{A}\cdot\mathbf{B}^{2}\right)$
and $\mathrm{tr}\left(\mathbf{A}^{2}\cdot\mathbf{B}^{2}\right)$.

\subsection{Common Differential Operations in Tensor Notation}

$\bullet$ Here we present the most common differential operations
as defined by tensor notation. These operations are mostly based on
the various types of interaction between the vector differential operator
nabla $\nabla$ with tensors of different ranks as well as interaction
with other types of operation like dot and cross products.

$\bullet$ $\nabla$ is essentially a spatial partial differential
operator defined in Cartesian coordinate systems by:
\begin{equation}
\nabla_{i}=\frac{\partial}{\partial x_{i}}
\end{equation}

The definition of $\nabla$ in some non-Cartesian systems will be
given in $\S$ \ref{subOtherCoordinateSystems}.

\subsubsection{Cartesian System}

$\bullet$ The gradient of a differentiable scalar function of position
$f$ is a vector given by:
\begin{equation}
\left[\nabla f\right]_{i}=\nabla_{i}f=\frac{\partial f}{\partial x_{i}}=\partial_{i}f=f_{,i}\label{eqGrad}
\end{equation}

$\bullet$ The gradient of a differentiable vector function of position
$\mathbf{A}$ (which is the outer product, as defined in $\S$ \ref{subTensorMultiplication},
between the $\nabla$ operator and the vector) is a rank-2 tensor
defined by:
\begin{equation}
\left[\nabla\mathbf{A}\right]_{ij}=\partial_{i}A_{j}\label{eqGrad2}
\end{equation}

$\bullet$ The gradient operation is distributive but not commutative
or associative:
\begin{equation}
\nabla\left(f+h\right)=\nabla f+\nabla h
\end{equation}
\begin{equation}
\nabla f\ne f\nabla
\end{equation}
\begin{equation}
\left(\nabla f\right)h\ne\nabla\left(fh\right)
\end{equation}
where $f$ and $h$ are differentiable scalar functions of position.

$\bullet$ The divergence of a differentiable vector $\mathbf{A}$
is a scalar given by:
\begin{equation}
\nabla\cdot\mathbf{A}=\delta_{ij}\frac{\partial A_{i}}{\partial x_{j}}=\frac{\partial A_{i}}{\partial x_{i}}=\nabla_{i}A_{i}=\partial_{i}A_{i}=A_{i,i}\label{eqDiv}
\end{equation}
The divergence operation can also be viewed as taking the gradient
of the vector followed by a contraction. Hence, the divergence of
a vector is invariant because it is the trace of a rank-2 tensor.\footnote{It may also be argued that the divergence of a vector is a scalar
and hence it is invariant.}

$\bullet$ The divergence of a differentiable rank-2 tensor $\mathbf{A}$
is a vector defined in one of its forms by:
\begin{equation}
\left[\nabla\cdot\mathbf{A}\right]_{i}=\partial_{j}A_{ji}\label{eqDiv2}
\end{equation}
and in another form by
\begin{equation}
\left[\nabla\cdot\mathbf{A}\right]_{j}=\partial_{i}A_{ji}\label{eqDiv3}
\end{equation}
These two different forms can be given, respectively, in symbolic
notation by:
\begin{equation}
\nabla\cdot\mathbf{A}\,\,\,\,\,\,\,\,\,\,\&\,\,\,\,\,\,\,\,\,\,\nabla\cdot\mathbf{A}^{T}
\end{equation}
where $\mathbf{A}^{T}$ is the transpose of $\mathbf{A}$. More generally,
the divergence of a tensor of rank $n\ge2$, which is a tensor of
rank-($n-1$), can be defined in several forms, which are different
in general, depending on the combination of the contracted indices.

$\bullet$ The divergence operation is distributive but not commutative
or associative:
\begin{equation}
\nabla\cdot\left(\mathbf{A}+\mathbf{B}\right)=\nabla\cdot\mathbf{A}+\nabla\cdot\mathbf{B}
\end{equation}
\begin{equation}
\nabla\cdot\mathbf{A}\ne\mathbf{A}\cdot\nabla
\end{equation}
\begin{equation}
\nabla\cdot\left(f\mathbf{A}\right)\ne\nabla f\cdot\mathbf{A}
\end{equation}
where $\mathbf{A}$ and $\mathbf{B}$ are differentiable tensor functions
of position.

$\bullet$ The curl of a differentiable vector $\mathbf{A}$ is a
vector given by:
\begin{equation}
\left[\nabla\times\mathbf{A}\right]_{i}=\epsilon_{ijk}\frac{\partial A_{k}}{\partial x_{j}}=\epsilon_{ijk}\nabla_{j}A_{k}=\epsilon_{ijk}\partial_{j}A_{k}=\epsilon_{ijk}A_{k,j}\label{EqCurl}
\end{equation}

$\bullet$ The curl operation may be generalized to tensors of rank
$>1$, and hence the curl of a differentiable rank-2 tensor $\mathbf{A}$
can be defined as a rank-2 tensor given by:
\begin{equation}
\left[\nabla\times\mathbf{A}\right]_{ij}=\epsilon_{imn}\partial_{m}A_{nj}
\end{equation}

$\bullet$ The curl operation is distributive but not commutative
or associative:
\begin{equation}
\nabla\times\left(\mathbf{A}+\mathbf{B}\right)=\nabla\times\mathbf{A}+\nabla\times\mathbf{B}
\end{equation}
\begin{equation}
\nabla\times\mathbf{A}\ne\mathbf{A}\times\nabla
\end{equation}
\begin{equation}
\nabla\times\left(\mathbf{A\times B}\right)\ne\left(\nabla\times\mathbf{A}\right)\times\mathbf{B}
\end{equation}

$\bullet$ The Laplacian scalar operator, also called the harmonic
operator, acting on a differentiable scalar $f$ is given by:
\begin{equation}
\Delta f=\nabla^{2}f=\delta_{ij}\frac{\partial^{2}f}{\partial x_{i}\partial x_{j}}=\frac{\partial^{2}f}{\partial x_{i}\partial x_{i}}=\nabla_{ii}f=\partial_{ii}f=f_{,ii}\label{eqLaplacian}
\end{equation}

$\bullet$ The Laplacian operator acting on a differentiable vector
$\mathbf{A}$ is defined for each component of the vector similar
to the definition of the Laplacian acting on a scalar, that is
\begin{equation}
\left[\nabla^{2}\mathbf{A}\right]_{i}=\partial_{jj}A_{i}\label{eqLaplacian2}
\end{equation}

$\bullet$ The following scalar differential operator is commonly
used in science (e.g. in fluid dynamics):
\begin{equation}
\mathbf{A}\cdot\nabla=A_{i}\nabla_{i}=A_{i}\frac{\partial}{\partial x_{i}}=A_{i}\partial_{i}\label{eqANabla}
\end{equation}
where $\mathbf{A}$ is a vector. As indicated earlier, the order of
$A_{i}$ and $\partial_{i}$ should be respected.

$\bullet$ The following vector differential operator also has common
applications in science:
\begin{equation}
\left[\mathbf{A}\times\nabla\right]_{i}=\epsilon_{ijk}A_{j}\partial_{k}
\end{equation}

$\bullet$ The differentiation of a tensor increases its rank by one,
by introducing an extra covariant index, unless it implies a contraction
in which case it reduces the rank by one. Therefore the gradient of
a scalar is a vector and the gradient of a vector is a rank-2 tensor
($\partial_{i}A_{j}$), while the divergence of a vector is a scalar
and the divergence of a rank-2 tensor is a vector ($\partial_{j}A_{ji}$
or $\partial_{i}A_{ji}$). This may be justified by the fact that
$\nabla$ is a vector operator. On the other hand the Laplacian operator
does not change the rank since it is a scalar operator; hence the
Laplacian of a scalar is a scalar and the Laplacian of a vector is
a vector.

\subsubsection{Other Coordinate Systems\label{subOtherCoordinateSystems}}

$\bullet$ For completeness, we define here some differential operations
in the most commonly used non-Cartesian coordinate systems, namely
cylindrical and spherical systems, as well as general orthogonal coordinate
systems.

$\bullet$ We can use indexed generalized coordinates like $q_{1}$,
$q_{2}$ and $q_{3}$ for the cylindrical coordinates ($\rho,\phi,z$)
and the spherical coordinates ($r,\theta,\phi$). However, for more
clarity at this level and to follow the more conventional practice,
we use the coordinates of these systems as suffixes in place of the
indices used in the tensor notation.\footnote{There is another reason that is these are physical components not
covariant or contravariant.}

$\bullet$ For the cylindrical system identified by the coordinates
($\rho,\phi,z$) with an orthonormal basis vectors $\mathbf{e}_{\rho},\,\mathbf{e}_{\phi}$
and $\mathbf{e}_{z}$:\footnote{It should be obvious that since $\rho,\,\phi$ and $z$ are specific
coordinates and not variable indices, the summation convention does
not apply. }

The $\nabla$ operator is:
\begin{equation}
\nabla=\mathbf{e}_{\rho}\partial_{\rho}+\mathbf{e}_{\phi}\frac{1}{\rho}\partial_{\phi}+\mathbf{e}_{z}\partial_{z}
\end{equation}
The Laplacian operator is:
\begin{equation}
\nabla^{2}=\partial_{\rho\rho}+\frac{1}{\rho}\partial_{\rho}+\frac{1}{\rho^{2}}\partial_{\phi\phi}+\partial_{zz}
\end{equation}
The gradient of a differentiable scalar $f$ is:
\begin{equation}
\nabla f=\mathbf{e}_{\rho}\frac{\partial f}{\partial\rho}+\mathbf{e}_{\phi}\frac{1}{\rho}\frac{\partial f}{\partial\phi}+\mathbf{e}_{z}\frac{\partial f}{\partial z}
\end{equation}
The divergence of a differentiable vector $\mathbf{A}$ is:
\begin{equation}
\nabla\cdot\mathbf{A}=\frac{1}{\rho}\left[\frac{\partial\left(\rho A_{\rho}\right)}{\partial\rho}+\frac{\partial A_{\phi}}{\partial\phi}+\frac{\partial\left(\rho A_{z}\right)}{\partial z}\right]
\end{equation}
The curl of a differentiable vector $\mathbf{A}$ is:
\begin{equation}
\nabla\times\mathbf{A}=\frac{1}{\rho}\begin{vmatrix}\begin{array}{ccc}
\mathbf{e}_{\rho} & \rho\mathbf{e}_{\phi} & \mathbf{e}_{z}\\
\frac{\partial}{\partial\rho} & \frac{\partial}{\partial\phi} & \frac{\partial}{\partial z}\\
A_{\rho} & \rho A_{\phi} & A_{z}
\end{array}\end{vmatrix}
\end{equation}
For plane polar coordinate systems, these operators and operations
can be obtained by dropping the $z$ components or terms from the
cylindrical form of the above operators and operations.

$\bullet$ For the spherical system identified by the coordinates
($r,\theta,\phi$) with an orthonormal basis vectors $\mathbf{e}_{r},\,\mathbf{e}_{\theta}$
and $\mathbf{e}_{\phi}$:\footnote{Again, the summation convention does not apply to $r,\,\theta$ and
$\phi$. }

The $\nabla$ operator is:
\begin{equation}
\nabla=\mathbf{e}_{r}\partial_{r}+\mathbf{e}_{\theta}\frac{1}{r}\partial_{\theta}+\mathbf{e}_{\phi}\frac{1}{r\sin\theta}\partial_{\phi}
\end{equation}
The Laplacian operator is:
\begin{equation}
\nabla^{2}=\partial_{rr}+\frac{2}{r}\partial_{r}+\frac{1}{r^{2}}\partial_{\theta\theta}+\frac{\cos\theta}{r^{2}\sin\theta}\partial_{\theta}+\frac{1}{r^{2}\sin^{2}\theta}\partial_{\phi\phi}
\end{equation}
The gradient of a differentiable scalar $f$ is:
\begin{equation}
\nabla f=\mathbf{e}_{r}\frac{\partial f}{\partial r}+\mathbf{e}_{\theta}\frac{1}{r}\frac{\partial f}{\partial\theta}+\mathbf{e}_{\phi}\frac{1}{r\sin\theta}\frac{\partial f}{\partial\phi}
\end{equation}
The divergence of a differentiable vector $\mathbf{A}$ is:
\begin{equation}
\nabla\cdot\mathbf{A}=\frac{1}{r^{2}\sin\theta}\left[\sin\theta\frac{\partial\left(r^{2}A_{r}\right)}{\partial r}+r\frac{\partial\left(\sin\theta A_{\theta}\right)}{\partial\theta}+r\frac{\partial A_{\phi}}{\partial\phi}\right]
\end{equation}
The curl of a differentiable vector $\mathbf{A}$ is:
\begin{equation}
\nabla\times\mathbf{A}=\frac{1}{r^{2}\sin\theta}\begin{vmatrix}\begin{array}{ccc}
\mathbf{e}_{r} & r\mathbf{e}_{\theta} & r\sin\theta\mathbf{e}_{\phi}\\
\frac{\partial}{\partial r} & \frac{\partial}{\partial\theta} & \frac{\partial}{\partial\phi}\\
A_{r} & rA_{\theta} & r\sin\theta A_{\phi}
\end{array}\end{vmatrix}
\end{equation}

$\bullet$ For a general orthogonal system in a 3D space identified
by the coordinates ($u_{1},u_{2},u_{3}$) with unit basis vectors
$\mathbf{u}_{1},\,\mathbf{u}_{2}$ and $\mathbf{u}_{3}$ and scale
factors $h_{1},\,h_{2}$ and $h_{3}$ where $h_{i}=\left|\frac{\partial\mathbf{r}}{\partial u_{i}}\right|$
and $\mathbf{r}$ is the position vector:

The $\nabla$ operator is:
\begin{equation}
\nabla=\frac{\mathbf{u}_{1}}{h_{1}}\frac{\partial}{\partial u_{1}}+\frac{\mathbf{u}_{2}}{h_{2}}\frac{\partial}{\partial u_{2}}+\frac{\mathbf{u}_{3}}{h_{3}}\frac{\partial}{\partial u_{3}}
\end{equation}
The Laplacian operator is:
\begin{equation}
\nabla^{2}=\frac{1}{h_{1}h_{2}h_{3}}\left[\frac{\partial}{\partial u_{1}}\left(\frac{h_{2}h_{3}}{h_{1}}\frac{\partial}{\partial u_{1}}\right)+\frac{\partial}{\partial u_{2}}\left(\frac{h_{1}h_{3}}{h_{2}}\frac{\partial}{\partial u_{2}}\right)+\frac{\partial}{\partial u_{3}}\left(\frac{h_{1}h_{2}}{h_{3}}\frac{\partial}{\partial u_{3}}\right)\right]
\end{equation}
The gradient of a differentiable scalar $f$ is:
\begin{equation}
\nabla f=\frac{\mathbf{u}_{1}}{h_{1}}\frac{\partial f}{\partial u_{1}}+\frac{\mathbf{u}_{2}}{h_{2}}\frac{\partial f}{\partial u_{2}}+\frac{\mathbf{u}_{3}}{h_{3}}\frac{\partial f}{\partial u_{3}}
\end{equation}
The divergence of a differentiable vector $\mathbf{A}$ is:
\begin{equation}
\nabla\cdot\mathbf{A}=\frac{1}{h_{1}h_{2}h_{3}}\left[\frac{\partial}{\partial u_{1}}\left(h_{2}h_{3}A_{1}\right)+\frac{\partial}{\partial u_{2}}\left(h_{1}h_{3}A_{2}\right)+\frac{\partial}{\partial u_{3}}\left(h_{1}h_{2}A_{3}\right)\right]
\end{equation}
The curl of a differentiable vector $\mathbf{A}$ is:
\begin{equation}
\nabla\times\mathbf{A}=\frac{1}{h_{1}h_{2}h_{3}}\begin{vmatrix}\begin{array}{ccc}
h_{1}\mathbf{u}_{1} & h_{2}\mathbf{u}_{2} & h_{3}\mathbf{u}_{3}\\
\frac{\partial}{\partial u_{1}} & \frac{\partial}{\partial u_{2}} & \frac{\partial}{\partial u_{3}}\\
h_{1}A_{1} & h_{2}A_{2} & h_{3}A_{3}
\end{array}\end{vmatrix}
\end{equation}

\subsection{Common Identities in Vector and Tensor Notation}

$\bullet$ Here we present some of the widely used identities of vector
calculus in the traditional vector notation and in its equivalent
tensor notation. In the following bullet points, $f$ and $h$ are
differentiable scalar fields; $\mathbf{A}$, $\mathbf{B}$, $\mathbf{C}$
and $\mathbf{D}$ are differentiable vector fields; and $\mathbf{r}=x_{i}\mathbf{e}_{i}$
is the position vector.

$\bullet$
\begin{eqnarray}
\nabla\cdot\mathbf{r} & = & n\nonumber \\
 & \Updownarrow\\
\partial_{i}x_{i} & = & n\nonumber
\end{eqnarray}
where $n$ is the space dimension.

$\bullet$
\begin{eqnarray}
\nabla\times\mathbf{r} & = & \mathbf{0}\nonumber \\
 & \Updownarrow\\
\epsilon_{ijk}\partial_{j}x_{k} & = & 0\nonumber
\end{eqnarray}

$\bullet$
\begin{eqnarray}
\nabla\left(\mathbf{a}\cdot\mathbf{r}\right) & = & \mathbf{a}\nonumber \\
 & \Updownarrow\\
\partial_{i}\left(a_{j}x_{j}\right) & = & a_{i}\nonumber
\end{eqnarray}
where $\mathbf{a}$ is a constant vector.

$\bullet$
\begin{eqnarray}
\nabla\cdot\left(\nabla f\right) & = & \nabla^{2}f\nonumber \\
 & \Updownarrow\\
\partial_{i}\left(\partial_{i}f\right) & = & \partial_{ii}f\nonumber
\end{eqnarray}

$\bullet$
\begin{eqnarray}
\nabla\cdot\left(\nabla\times\mathbf{A}\right) & = & 0\nonumber \\
 & \Updownarrow\\
\epsilon_{ijk}\partial_{i}\partial_{j}A_{k} & = & 0\nonumber
\end{eqnarray}

$\bullet$
\begin{eqnarray}
\nabla\times\left(\nabla f\right) & = & \mathbf{0}\nonumber \\
 & \Updownarrow\\
\epsilon_{ijk}\partial_{j}\partial_{k}f & = & 0\nonumber
\end{eqnarray}

$\bullet$
\begin{eqnarray}
\nabla\left(fh\right) & = & f\nabla h+h\nabla f\nonumber \\
 & \Updownarrow\\
\partial_{i}\left(fh\right) & = & f\partial_{i}h+h\partial_{i}f\nonumber
\end{eqnarray}

$\bullet$
\begin{eqnarray}
\nabla\cdot\left(f\mathbf{A}\right) & = & f\nabla\cdot\mathbf{A}+\mathbf{A}\cdot\nabla f\nonumber \\
 & \Updownarrow\\
\partial_{i}\left(fA_{i}\right) & = & f\partial_{i}A_{i}+A_{i}\partial_{i}f\nonumber
\end{eqnarray}

$\bullet$
\begin{eqnarray}
\nabla\times\left(f\mathbf{A}\right) & = & f\nabla\times\mathbf{A}+\nabla f\times\mathbf{A}\nonumber \\
 & \Updownarrow\\
\epsilon_{ijk}\partial_{j}\left(fA_{k}\right) & = & f\epsilon_{ijk}\partial_{j}A_{k}+\epsilon_{ijk}\left(\partial_{j}f\right)A_{k}\nonumber
\end{eqnarray}

$\bullet$
\begin{alignat}{3}
\mathbf{A}\cdot\left(\mathbf{B}\times\mathbf{C}\right) & = & \,\mathbf{C}\cdot\left(\mathbf{A}\times\mathbf{B}\right) & = & \,\mathbf{B}\cdot\left(\mathbf{C}\times\mathbf{A}\right)\nonumber \\
 & \Updownarrow &  & \Updownarrow\\
\epsilon_{ijk}A_{i}B_{j}C_{k} & = & \epsilon_{kij}C_{k}A_{i}B_{j} & = & \epsilon_{jki}B_{j}C_{k}A_{i}\nonumber
\end{alignat}

$\bullet$
\begin{eqnarray}
\mathbf{A}\times\left(\mathbf{B}\times\mathbf{C}\right) & = & \mathbf{B}\left(\mathbf{A}\cdot\mathbf{C}\right)-\mathbf{C}\left(\mathbf{A}\cdot\mathbf{B}\right)\nonumber \\
 & \Updownarrow\\
\epsilon_{ijk}A_{j}\epsilon_{klm}B_{l}C_{m} & = & B_{i}\left(A_{m}C_{m}\right)-C_{i}\left(A_{l}B_{l}\right)\nonumber
\end{eqnarray}

$\bullet$
\begin{eqnarray}
\mathbf{A}\times\left(\nabla\times\mathbf{B}\right) & = & \left(\nabla\mathbf{B}\right)\cdot\mathbf{A}-\mathbf{A}\cdot\nabla\mathbf{B}\nonumber \\
 & \Updownarrow\\
\epsilon_{ijk}\epsilon_{klm}A_{j}\partial_{l}B_{m} & = & \left(\partial_{i}B_{m}\right)A_{m}-A_{l}\left(\partial_{l}B_{i}\right)\nonumber
\end{eqnarray}

$\bullet$
\begin{eqnarray}
\nabla\times\left(\nabla\times\mathbf{A}\right) & = & \nabla\left(\nabla\cdot\mathbf{A}\right)-\nabla^{2}\mathbf{A}\nonumber \\
 & \Updownarrow\\
\epsilon_{ijk}\epsilon_{klm}\partial_{j}\partial_{l}A_{m} & = & \partial_{i}\left(\partial_{m}A_{m}\right)-\partial_{ll}A_{i}\nonumber
\end{eqnarray}

$\bullet$
\begin{eqnarray}
\nabla\left(\mathbf{A}\cdot\mathbf{B}\right) & = & \mathbf{A}\times\left(\nabla\times\mathbf{B}\right)+\mathbf{B}\times\left(\nabla\times\mathbf{A}\right)+\left(\mathbf{A}\cdot\nabla\right)\mathbf{B}+\left(\mathbf{B}\cdot\nabla\right)\mathbf{A}\nonumber \\
 & \Updownarrow\\
\partial_{i}\left(A_{m}B_{m}\right) & = & \epsilon_{ijk}A_{j}\left(\epsilon_{klm}\partial_{l}B_{m}\right)+\epsilon_{ijk}B_{j}\left(\epsilon_{klm}\partial_{l}A_{m}\right)+\left(A_{l}\partial_{l}\right)B_{i}+\left(B_{l}\partial_{l}\right)A_{i}\nonumber
\end{eqnarray}

$\bullet$
\begin{eqnarray}
\nabla\cdot\left(\mathbf{A}\times\mathbf{B}\right) & = & \mathbf{B}\cdot\left(\nabla\times\mathbf{A}\right)-\mathbf{A}\cdot\left(\nabla\times\mathbf{B}\right)\nonumber \\
 & \Updownarrow\\
\partial_{i}\left(\epsilon_{ijk}A_{j}B_{k}\right) & = & B_{k}\left(\epsilon_{kij}\partial_{i}A_{j}\right)-A_{j}\left(\epsilon_{jik}\partial_{i}B_{k}\right)\nonumber
\end{eqnarray}

$\bullet$
\begin{eqnarray}
\nabla\times\left(\mathbf{A}\times\mathbf{B}\right) & = & \left(\mathbf{B}\cdot\nabla\right)\mathbf{A}+\left(\nabla\cdot\mathbf{B}\right)\mathbf{A}-\left(\nabla\cdot\mathbf{A}\right)\mathbf{B}-\left(\mathbf{A}\cdot\nabla\right)\mathbf{B}\nonumber \\
 & \Updownarrow\\
\epsilon_{ijk}\epsilon_{klm}\partial_{j}\left(A_{l}B_{m}\right) & = & \left(B_{m}\partial_{m}\right)A_{i}+\left(\partial_{m}B_{m}\right)A_{i}-\left(\partial_{j}A_{j}\right)B_{i}-\left(A_{j}\partial_{j}\right)B_{i}\nonumber
\end{eqnarray}

$\bullet$
\begin{eqnarray}
\left(\mathbf{A}\times\mathbf{B}\right)\cdot\left(\mathbf{C}\times\mathbf{D}\right) & = & \begin{vmatrix}\begin{array}{cc}
\mathbf{A}\cdot\mathbf{C} & \mathbf{A}\cdot\mathbf{D}\\
\mathbf{B}\cdot\mathbf{C} & \mathbf{B}\cdot\mathbf{D}
\end{array}\end{vmatrix}\nonumber \\
 & \Updownarrow\\
\epsilon_{ijk}A_{j}B_{k}\epsilon_{ilm}C_{l}D_{m} & = & \left(A_{l}C_{l}\right)\left(B_{m}D_{m}\right)-\left(A_{m}D_{m}\right)\left(B_{l}C_{l}\right)\nonumber
\end{eqnarray}

$\bullet$
\begin{eqnarray}
\left(\mathbf{A}\times\mathbf{B}\right)\times\left(\mathbf{C}\times\mathbf{D}\right) & = & \left[\mathbf{D}\cdot\left(\mathbf{A}\times\mathbf{B}\right)\right]\mathbf{C}-\left[\mathbf{C}\cdot\left(\mathbf{A}\times\mathbf{B}\right)\right]\mathbf{D}\nonumber \\
 & \Updownarrow\\
\epsilon_{ijk}\epsilon_{jmn}A_{m}B_{n}\epsilon_{kpq}C_{p}D_{q} & = & \left(\epsilon_{qmn}D_{q}A_{m}B_{n}\right)C_{i}-\left(\epsilon_{pmn}C_{p}A_{m}B_{n}\right)D_{i}\nonumber
\end{eqnarray}

$\bullet$ In vector and tensor notations, the condition for a vector
field $\mathbf{A}$ to be solenoidal is:
\begin{eqnarray}
\nabla\cdot\mathbf{A} & = & 0\nonumber \\
 & \Updownarrow\\
\partial_{i}A_{i} & = & 0\nonumber
\end{eqnarray}

$\bullet$ In vector and tensor notations, the condition for a vector
field $\mathbf{A}$ to be irrotational is:
\begin{eqnarray}
\nabla\times\mathbf{A} & = & \mathbf{0}\nonumber \\
 & \Updownarrow\\
\epsilon_{ijk}\partial_{j}A_{k} & = & 0\nonumber
\end{eqnarray}

\subsection{Integral Theorems in Tensor Notation}

$\bullet$ The divergence theorem for a differentiable vector field
$\mathbf{A}$ in vector and tensor notation is:
\begin{eqnarray}
\iiint_{V}\nabla\cdot\mathbf{A}\,d\tau & = & \iint_{S}\mathbf{A}\cdot\mathbf{n}\,d\sigma\nonumber \\
 & \Updownarrow\\
\int_{V}\partial_{i}A_{i}d\tau & = & \int_{S}A_{i}n_{i}d\sigma\nonumber
\end{eqnarray}
where $V$ is a bounded region in an $n$D space enclosed by a generalized
surface $S$, $d\tau$ and $d\sigma$ are generalized volume and surface
elements respectively, $\mathbf{n}$ and $n_{i}$ are unit normal
to the surface and its $i^{th}$ component respectively, and the index
$i$ ranges over $1,\ldots,n$.

$\bullet$ The divergence theorem for a differentiable rank-2 tensor
field $\mathbf{A}$ in tensor notation for the first index is given
by:
\begin{equation}
\int_{V}\partial_{i}A_{il}d\tau=\int_{S}A_{il}n_{i}d\sigma
\end{equation}

$\bullet$ The divergence theorem for differentiable tensor fields
of higher ranks $\mathbf{A}$ in tensor notation for the index $k$
is:
\begin{equation}
\int_{V}\partial_{k}A_{ij\ldots k\ldots m}d\tau=\int_{S}A_{ij\ldots k\ldots m}n_{k}d\sigma
\end{equation}

$\bullet$ Stokes theorem for a differentiable vector field $\mathbf{A}$
in vector and tensor notation is:
\begin{eqnarray}
\iint_{S}\left(\nabla\times\mathbf{A}\right)\cdot\mathbf{n}\,d\sigma & = & \int_{C}\mathbf{A}\cdot d\mathbf{r}\nonumber \\
 & \Updownarrow\\
\int_{S}\epsilon_{ijk}\partial_{j}A_{k}n_{i}d\sigma & = & \int_{C}A_{i}dx_{i}\nonumber
\end{eqnarray}
where $C$ stands for the perimeter of the surface $S$ and $d\mathbf{r}$
is the vector element tangent to the perimeter.

$\bullet$ Stokes theorem for a differentiable rank-2 tensor field
$\mathbf{A}$ in tensor notation for the first index is:
\begin{equation}
\int_{S}\epsilon_{ijk}\partial_{j}A_{kl}n_{i}d\sigma=\int_{C}A_{il}dx_{i}
\end{equation}

$\bullet$ Stokes theorem for differentiable tensor fields of higher
ranks $\mathbf{A}$ in tensor notation for the index $k$ is:
\begin{equation}
\int_{S}\epsilon_{ijk}\partial_{j}A_{lm\ldots k\ldots n}n_{i}d\sigma=\int_{C}A_{lm\ldots k\ldots n}dx_{k}
\end{equation}

\subsection{Examples of Using Tensor Techniques to Prove Identities\label{subProvingIdentities}}

$\bullet$ $\nabla\cdot\mathbf{r}=n$:
\begin{equation}
\begin{aligned}\nabla\cdot\mathbf{r} & =\partial_{i}x_{i} & \,\,\,\,\,\,\,\,\,\,\,\,\,\, & \text{(Eq. \ref{eqDiv})}\\
 & =\delta_{ii} &  & \text{(Eq. \ref{eqdxn})}\\
 & =n &  & \text{(Eq. \ref{eqdxn})}
\end{aligned}
\end{equation}

$\bullet$ $\nabla\times\mathbf{r}=\mathbf{0}$:
\begin{equation}
\begin{aligned}\left[\nabla\times\mathbf{r}\right]_{i} & =\epsilon_{ijk}\partial_{j}x_{k} & \,\,\,\,\,\,\,\,\,\,\,\,\,\,\,\,\, & \text{(Eq. \ref{EqCurl})}\\
 & =\epsilon_{ijk}\delta_{kj} &  & \text{(Eq. \ref{eqdxdelta})}\\
 & =\epsilon_{ijj} &  & \text{(Eq. \ref{EqIndexReplace})}\\
 & =0 &  & \text{(Eq. \ref{eqEpsilon3Definition})}
\end{aligned}
\end{equation}
Since $i$ is a free index the identity is proved for all components.

$\bullet$ $\nabla\left(\mathbf{a}\cdot\mathbf{r}\right)=\mathbf{a}$:
\begin{equation}
\begin{aligned}\left[\nabla\left(\mathbf{a}\cdot\mathbf{r}\right)\right]_{i} & =\partial_{i}\left(a_{j}x_{j}\right) & \,\,\,\,\,\,\,\,\,\,\,\,\,\,\,\,\,\,\,\, & \text{(Eqs. \ref{eqGrad} \& \ref{eqDotProduct})}\\
 & =a_{j}\partial_{i}x_{j}+x_{j}\partial_{i}a_{j} &  & \text{(product rule)}\\
 & =a_{j}\partial_{i}x_{j} &  & \text{(\ensuremath{a_{j}} is constant)}\\
 & =a_{j}\delta_{ji} &  & \text{(Eq. \ref{eqdxdelta})}\\
 & =a_{i} &  & \text{(Eq. \ref{EqIndexReplace})}\\
 & =\left[\mathbf{a}\right]_{i} &  & \text{(definition of index)}
\end{aligned}
\end{equation}
Since $i$ is a free index the identity is proved for all components.

$\bullet$ $\nabla\cdot\left(\nabla f\right)=\nabla^{2}f$:
\begin{equation}
\begin{aligned}\nabla\cdot\left(\nabla f\right) & =\partial_{i}\left[\nabla f\right]_{i} & \,\,\,\,\,\,\,\,\,\,\,\,\,\,\,\,\,\, & \text{(Eq. \ref{eqDiv})}\\
 & =\partial_{i}\left(\partial_{i}f\right) &  & \text{(Eq. \ref{eqGrad})}\\
 & =\partial_{i}\partial_{i}f &  & \text{(rules of differentiation)}\\
 & =\partial_{ii}f &  & \text{(definition of 2nd derivative)}\\
 & =\nabla^{2}f &  & \text{(Eq. \ref{eqLaplacian})}
\end{aligned}
\end{equation}

$\bullet$ $\nabla\cdot\left(\nabla\times\mathbf{A}\right)=0$:
\begin{equation}
\begin{aligned}\nabla\cdot\left(\nabla\times\mathbf{A}\right) & =\partial_{i}\left[\nabla\times\mathbf{A}\right]_{i} & \,\,\,\,\,\,\,\,\,\,\,\,\,\,\,\, & \text{(Eq. \ref{eqDiv})}\\
 & =\partial_{i}\left(\epsilon_{ijk}\partial_{j}A_{k}\right) &  & \text{(Eq. \ref{EqCurl})}\\
 & =\epsilon_{ijk}\partial_{i}\partial_{j}A_{k} &  & \text{(\ensuremath{\partial}\ not acting on \ensuremath{\epsilon})}\\
 & =\epsilon_{ijk}\partial_{j}\partial_{i}A_{k} &  & \text{(continuity condition)}\\
 & =-\epsilon_{jik}\partial_{j}\partial_{i}A_{k} &  & \text{(Eq. \ref{EqEpsilonCycle})}\\
 & =-\epsilon_{ijk}\partial_{i}\partial_{j}A_{k} &  & \text{(relabeling dummy indices \ensuremath{i} and \ensuremath{j})}\\
 & =0 &  & \text{(since \ensuremath{\epsilon_{ijk}\partial_{i}\partial_{j}A_{k}=-\epsilon_{ijk}\partial_{i}\partial_{j}A_{k}})}
\end{aligned}
\end{equation}
This can also be concluded from line three by arguing that: since
by the continuity condition $\partial_{i}$ and $\partial_{j}$ can
change their order with no change in the value of the term while a
corresponding change of the order of $i$ and $j$ in $\epsilon_{ijk}$
results in a sign change, we see that each term in the sum has its
own negative and hence the terms add up to zero (see Eq. \ref{eqPermutingTwoFactors}).

$\bullet$ $\nabla\times\left(\nabla f\right)=\mathbf{0}$:
\begin{equation}
\begin{aligned}\left[\nabla\times\left(\nabla f\right)\right]_{i} & =\epsilon_{ijk}\partial_{j}\left[\nabla f\right]_{k} & \,\,\,\,\,\,\,\,\,\,\,\,\,\,\,\,\,\,\, & \text{(Eq. \ref{EqCurl})}\\
 & =\epsilon_{ijk}\partial_{j}\left(\partial_{k}f\right) &  & \text{(Eq. \ref{eqGrad})}\\
 & =\epsilon_{ijk}\partial_{j}\partial_{k}f &  & \text{(rules of differentiation)}\\
 & =\epsilon_{ijk}\partial_{k}\partial_{j}f &  & \text{(continuity condition)}\\
 & =-\epsilon_{ikj}\partial_{k}\partial_{j}f &  & \text{(Eq. \ref{EqEpsilonCycle})}\\
 & =-\epsilon_{ijk}\partial_{j}\partial_{k}f &  & \text{(relabeling dummy indices \ensuremath{j} and \ensuremath{k})}\\
 & =0 &  & \text{(since \ensuremath{\epsilon_{ijk}\partial_{j}\partial_{k}f=-\epsilon_{ijk}\partial_{j}\partial_{k}f})}
\end{aligned}
\end{equation}
This can also be concluded from line three by a similar argument to
the one given in the previous point. Because $\left[\nabla\times\left(\nabla f\right)\right]_{i}$
is an arbitrary component, then each component is zero.

$\bullet$ $\nabla\left(fh\right)=f\nabla h+h\nabla f$:
\begin{equation}
\begin{aligned}\left[\nabla\left(fh\right)\right]_{i} & =\partial_{i}\left(fh\right) & \,\,\,\,\,\,\,\,\,\,\,\,\,\,\,\,\, & \text{(Eq. \ref{eqGrad})}\\
 & =f\partial_{i}h+h\partial_{i}f &  & \text{(product rule)}\\
 & =\left[f\nabla h\right]_{i}+\left[h\nabla f\right]_{i} &  & \text{(Eq. \ref{eqGrad})}\\
 & =\left[f\nabla h+h\nabla f\right]_{i} &  & \text{(Eq. \ref{eqIndexDistributive1})}
\end{aligned}
\end{equation}
Because $i$ is a free index the identity is proved for all components.

$\bullet$ $\nabla\cdot\left(f\mathbf{A}\right)=f\nabla\cdot\mathbf{A}+\mathbf{A}\cdot\nabla f$:
\begin{equation}
\begin{aligned}\nabla\cdot\left(f\mathbf{A}\right) & =\partial_{i}\left[f\mathbf{A}\right]_{i} & \,\,\,\,\,\,\,\,\,\,\,\,\,\,\,\,\, & \text{(Eq. \ref{eqDiv})}\\
 & =\partial_{i}\left(fA_{i}\right) &  & \text{(definition of index)}\\
 & =f\partial_{i}A_{i}+A_{i}\partial_{i}f &  & \text{(product rule)}\\
 & =f\nabla\cdot\mathbf{A}+\mathbf{A}\cdot\nabla f &  & \text{(Eqs. \ref{eqDiv} \& \ref{eqANabla})}
\end{aligned}
\end{equation}

$\bullet$ $\nabla\times\left(f\mathbf{A}\right)=f\nabla\times\mathbf{A}+\nabla f\times\mathbf{A}$:
\begin{equation}
\begin{aligned}\left[\nabla\times\left(f\mathbf{A}\right)\right]_{i} & =\epsilon_{ijk}\partial_{j}\left[f\mathbf{A}\right]_{k} & \,\,\,\,\,\,\,\,\,\,\,\,\,\,\,\,\,\,\,\, & \text{(Eq. \ref{EqCurl})}\\
 & =\epsilon_{ijk}\partial_{j}\left(fA_{k}\right) &  & \text{(definition of index)}\\
 & =f\epsilon_{ijk}\partial_{j}A_{k}+\epsilon_{ijk}\left(\partial_{j}f\right)A_{k} &  & \text{(product rule \& commutativity)}\\
 & =f\epsilon_{ijk}\partial_{j}A_{k}+\epsilon_{ijk}\left[\nabla f\right]_{j}A_{k} &  & \text{(Eq. \ref{eqGrad})}\\
 & =\left[f\nabla\times\mathbf{A}\right]_{i}+\left[\nabla f\times\mathbf{A}\right]_{i} &  & \text{(Eqs. \ref{EqCurl} \& \ref{EqCrossProduct})}\\
 & =\left[f\nabla\times\mathbf{A}+\nabla f\times\mathbf{A}\right]_{i} &  & \text{(Eq. \ref{eqIndexDistributive1})}
\end{aligned}
\end{equation}
Because $i$ is a free index the identity is proved for all components.

$\bullet$ $\mathbf{A}\cdot\left(\mathbf{B}\times\mathbf{C}\right)=\mathbf{C}\cdot\left(\mathbf{A}\times\mathbf{B}\right)=\mathbf{B}\cdot\left(\mathbf{C}\times\mathbf{A}\right)$:
\begin{equation}
\begin{aligned}\mathbf{A}\cdot\left(\mathbf{B}\times\mathbf{C}\right) & =\epsilon_{ijk}A_{i}B_{j}C_{k} & \,\,\,\,\,\,\,\,\,\,\,\,\,\,\,\, & \text{(Eq. \ref{EqScalarTripleProduct})}\\
 & =\epsilon_{kij}A_{i}B_{j}C_{k} &  & \text{(Eq. \ref{EqEpsilonCycle})}\\
 & =\epsilon_{kij}C_{k}A_{i}B_{j} &  & \text{(commutativity)}\\
 & =\mathbf{C}\cdot\left(\mathbf{A}\times\mathbf{B}\right) &  & \text{(Eq. \ref{EqScalarTripleProduct})}\\
 & =\epsilon_{jki}A_{i}B_{j}C_{k} &  & \text{(Eq. \ref{EqEpsilonCycle})}\\
 & =\epsilon_{jki}B_{j}C_{k}A_{i} &  & \text{(commutativity)}\\
 & =\mathbf{B}\cdot\left(\mathbf{C}\times\mathbf{A}\right) &  & \text{(Eq. \ref{EqScalarTripleProduct})}
\end{aligned}
\end{equation}
The negative permutations of these identities can be similarly obtained
and proved by changing the order of the vectors in the cross products
which results in a sign change.

$\bullet$ $\mathbf{A}\times\left(\mathbf{B}\times\mathbf{C}\right)=\mathbf{B}\left(\mathbf{A}\cdot\mathbf{C}\right)-\mathbf{C}\left(\mathbf{A}\cdot\mathbf{B}\right)$:
\begin{equation}
\begin{aligned}\left[\mathbf{A}\times\left(\mathbf{B}\times\mathbf{C}\right)\right]_{i} & =\epsilon_{ijk}A_{j}\left[\mathbf{B}\times\mathbf{C}\right]_{k} & \,\,\,\,\,\, & \text{(Eq. \ref{EqCrossProduct})}\\
 & =\epsilon_{ijk}A_{j}\epsilon_{klm}B_{l}C_{m} &  & \text{(Eq. \ref{EqCrossProduct})}\\
 & =\epsilon_{ijk}\epsilon_{klm}A_{j}B_{l}C_{m} &  & \text{(commutativity)}\\
 & =\epsilon_{ijk}\epsilon_{lmk}A_{j}B_{l}C_{m} &  & \text{(Eq. \ref{EqEpsilonCycle})}\\
 & =\left(\delta_{il}\delta_{jm}-\delta_{im}\delta_{jl}\right)A_{j}B_{l}C_{m} &  & \text{(Eq. \ref{EqEpsilonDelta})}\\
 & =\delta_{il}\delta_{jm}A_{j}B_{l}C_{m}-\delta_{im}\delta_{jl}A_{j}B_{l}C_{m} &  & \text{(distributivity)}\\
 & =\left(\delta_{il}B_{l}\right)\left(\delta_{jm}A_{j}C_{m}\right)-\left(\delta_{im}C_{m}\right)\left(\delta_{jl}A_{j}B_{l}\right) &  & \text{(commutativity and grouping)}\\
 & =B_{i}\left(A_{m}C_{m}\right)-C_{i}\left(A_{l}B_{l}\right) &  & \text{(Eq. \ref{EqIndexReplace})}\\
 & =B_{i}\left(\mathbf{A}\cdot\mathbf{C}\right)-C_{i}\left(\mathbf{A}\cdot\mathbf{B}\right) &  & \text{(Eq. \ref{eqDotProduct})}\\
 & =\left[\mathbf{B}\left(\mathbf{A}\cdot\mathbf{C}\right)\right]_{i}-\left[\mathbf{C}\left(\mathbf{A}\cdot\mathbf{B}\right)\right]_{i} &  & \text{(definition of index)}\\
 & =\left[\mathbf{B}\left(\mathbf{A}\cdot\mathbf{C}\right)-\mathbf{C}\left(\mathbf{A}\cdot\mathbf{B}\right)\right]_{i} &  & \text{(Eq. \ref{eqIndexDistributive1})}
\end{aligned}
\end{equation}
Because $i$ is a free index the identity is proved for all components.
Other variants of this identity {[}e.g. $\left(\mathbf{A}\times\mathbf{B}\right)\times\mathbf{C}${]}
can be obtained and proved similarly by changing the order of the
factors in the external cross product with adding a minus sign.

$\bullet$ $\mathbf{A}\times\left(\nabla\times\mathbf{B}\right)=\left(\nabla\mathbf{B}\right)\cdot\mathbf{A}-\mathbf{A}\cdot\nabla\mathbf{B}$:
\begin{equation}
\begin{aligned}\left[\mathbf{A}\times\left(\nabla\times\mathbf{B}\right)\right]_{i} & =\epsilon_{ijk}A_{j}\left[\nabla\times\mathbf{B}\right]_{k} & \,\,\,\,\,\,\,\,\,\,\, & \text{(Eq. \ref{EqCrossProduct})}\\
 & =\epsilon_{ijk}A_{j}\epsilon_{klm}\partial_{l}B_{m} &  & \text{(Eq. \ref{EqCurl})}\\
 & =\epsilon_{ijk}\epsilon_{klm}A_{j}\partial_{l}B_{m} &  & \text{(commutativity)}\\
 & =\epsilon_{ijk}\epsilon_{lmk}A_{j}\partial_{l}B_{m} &  & \text{(Eq. \ref{EqEpsilonCycle})}\\
 & =\left(\delta_{il}\delta_{jm}-\delta_{im}\delta_{jl}\right)A_{j}\partial_{l}B_{m} &  & \text{(Eq. \ref{EqEpsilonDelta})}\\
 & =\delta_{il}\delta_{jm}A_{j}\partial_{l}B_{m}-\delta_{im}\delta_{jl}A_{j}\partial_{l}B_{m} &  & \text{(distributivity)}\\
 & =A_{m}\partial_{i}B_{m}-A_{l}\partial_{l}B_{i} &  & \text{(Eq. \ref{EqIndexReplace})}\\
 & =\left(\partial_{i}B_{m}\right)A_{m}-A_{l}\left(\partial_{l}B_{i}\right) &  & \text{(commutativity \& grouping)}\\
 & =\left[\left(\nabla\mathbf{B}\right)\cdot\mathbf{A}\right]_{i}-\left[\mathbf{A}\cdot\nabla\mathbf{B}\right]_{i} &  & \text{(Eq. \ref{eqGrad2} \& \S\ \ref{secInnerProduct})}\\
 & =\left[\left(\nabla\mathbf{B}\right)\cdot\mathbf{A}-\mathbf{A}\cdot\nabla\mathbf{B}\right]_{i} &  & \text{(Eq. \ref{eqIndexDistributive1})}
\end{aligned}
\end{equation}
Because $i$ is a free index the identity is proved for all components.

$\bullet$ $\nabla\times\left(\nabla\times\mathbf{A}\right)=\nabla\left(\nabla\cdot\mathbf{A}\right)-\nabla^{2}\mathbf{A}$:
\begin{equation}
\begin{aligned}\left[\nabla\times\left(\nabla\times\mathbf{A}\right)\right]_{i} & =\epsilon_{ijk}\partial_{j}\left[\nabla\times\mathbf{A}\right]_{k} & \,\,\,\,\,\,\,\,\,\,\, & \text{(Eq. \ref{EqCurl})}\\
 & =\epsilon_{ijk}\partial_{j}\left(\epsilon_{klm}\partial_{l}A_{m}\right) &  & \text{(Eq. \ref{EqCurl})}\\
 & =\epsilon_{ijk}\epsilon_{klm}\partial_{j}\left(\partial_{l}A_{m}\right) &  & \text{(\ensuremath{\partial} not acting on \ensuremath{\epsilon})}\\
 & =\epsilon_{ijk}\epsilon_{lmk}\partial_{j}\partial_{l}A_{m} &  & \text{(Eq. \ref{EqEpsilonCycle} \& definition of derivative)}\\
 & =\left(\delta_{il}\delta_{jm}-\delta_{im}\delta_{jl}\right)\partial_{j}\partial_{l}A_{m} &  & \text{(Eq. \ref{EqEpsilonDelta})}\\
 & =\delta_{il}\delta_{jm}\partial_{j}\partial_{l}A_{m}-\delta_{im}\delta_{jl}\partial_{j}\partial_{l}A_{m} &  & \text{(distributivity)}\\
 & =\partial_{m}\partial_{i}A_{m}-\partial_{l}\partial_{l}A_{i} &  & \text{(Eq. \ref{EqIndexReplace})}\\
 & =\partial_{i}\left(\partial_{m}A_{m}\right)-\partial_{ll}A_{i} &  & \text{(\ensuremath{\partial}\ shift, grouping \& Eq. \ref{eqLaplacianSymbol})}\\
 & =\left[\nabla\left(\nabla\cdot\mathbf{A}\right)\right]_{i}-\left[\nabla^{2}\mathbf{A}\right]_{i} &  & \text{(Eqs. \ref{eqDiv}, \ref{eqGrad} \& \ref{eqLaplacian2})}\\
 & =\left[\nabla\left(\nabla\cdot\mathbf{A}\right)-\nabla^{2}\mathbf{A}\right]_{i} &  & \text{(Eqs. \ref{eqIndexDistributive1})}
\end{aligned}
\end{equation}
Because $i$ is a free index the identity is proved for all components.
This identity can also be considered as an instance of the identity
before the last one, observing that in the second term on the right
hand side the Laplacian should precede the vector, and hence no independent
proof is required.

$\bullet$ $\nabla\left(\mathbf{A}\cdot\mathbf{B}\right)=\mathbf{A}\times\left(\nabla\times\mathbf{B}\right)+\mathbf{B}\times\left(\nabla\times\mathbf{A}\right)+\left(\mathbf{A}\cdot\nabla\right)\mathbf{B}+\left(\mathbf{B}\cdot\nabla\right)\mathbf{A}$:

We start from the right hand side and end with the left hand side{\footnotesize{}
\begin{eqnarray}
\left[\mathbf{A}\times\left(\nabla\times\mathbf{B}\right)+\mathbf{B}\times\left(\nabla\times\mathbf{A}\right)+\left(\mathbf{A}\cdot\nabla\right)\mathbf{B}+\left(\mathbf{B}\cdot\nabla\right)\mathbf{A}\right]_{i} & =\nonumber \\
\left[\mathbf{A}\times\left(\nabla\times\mathbf{B}\right)\right]_{i}+\left[\mathbf{B}\times\left(\nabla\times\mathbf{A}\right)\right]_{i}+\left[\left(\mathbf{A}\cdot\nabla\right)\mathbf{B}\right]_{i}+\left[\left(\mathbf{B}\cdot\nabla\right)\mathbf{A}\right]_{i} & = & \,\,\text{(Eq. \ref{eqIndexDistributive1})}\nonumber \\
\epsilon_{ijk}A_{j}\left[\nabla\times\mathbf{B}\right]_{k}+\epsilon_{ijk}B_{j}\left[\nabla\times\mathbf{A}\right]_{k}+\left(A_{l}\partial_{l}\right)B_{i}+\left(B_{l}\partial_{l}\right)A_{i} & = & \,\,\text{(Eqs. \ref{EqCrossProduct}, \ref{eqDiv} \& indexing)}\nonumber \\
\epsilon_{ijk}A_{j}\left(\epsilon_{klm}\partial_{l}B_{m}\right)+\epsilon_{ijk}B_{j}\left(\epsilon_{klm}\partial_{l}A_{m}\right)+\left(A_{l}\partial_{l}\right)B_{i}+\left(B_{l}\partial_{l}\right)A_{i} & = & \,\,\text{(Eq. \ref{EqCurl})}\nonumber \\
\epsilon_{ijk}\epsilon_{klm}A_{j}\partial_{l}B_{m}+\epsilon_{ijk}\epsilon_{klm}B_{j}\partial_{l}A_{m}+\left(A_{l}\partial_{l}\right)B_{i}+\left(B_{l}\partial_{l}\right)A_{i} & = & \,\,\text{(commutativity)}\nonumber \\
\epsilon_{ijk}\epsilon_{lmk}A_{j}\partial_{l}B_{m}+\epsilon_{ijk}\epsilon_{lmk}B_{j}\partial_{l}A_{m}+\left(A_{l}\partial_{l}\right)B_{i}+\left(B_{l}\partial_{l}\right)A_{i} & = & \,\,\text{(Eq. \ref{EqEpsilonCycle})}\nonumber \\
\left(\delta_{il}\delta_{jm}-\delta_{im}\delta_{jl}\right)A_{j}\partial_{l}B_{m}+\left(\delta_{il}\delta_{jm}-\delta_{im}\delta_{jl}\right)B_{j}\partial_{l}A_{m}+\left(A_{l}\partial_{l}\right)B_{i}+\left(B_{l}\partial_{l}\right)A_{i} & = & \,\,\text{(Eq. \ref{EqEpsilonDelta})}\\
\hspace*{-1.5cm}\left(\delta_{il}\delta_{jm}A_{j}\partial_{l}B_{m}-\delta_{im}\delta_{jl}A_{j}\partial_{l}B_{m}\right)+\left(\delta_{il}\delta_{jm}B_{j}\partial_{l}A_{m}-\delta_{im}\delta_{jl}B_{j}\partial_{l}A_{m}\right)+\left(A_{l}\partial_{l}\right)B_{i}+\left(B_{l}\partial_{l}\right)A_{i} & = & \,\,\text{(distributivity)}\nonumber \\
\delta_{il}\delta_{jm}A_{j}\partial_{l}B_{m}-A_{l}\partial_{l}B_{i}+\delta_{il}\delta_{jm}B_{j}\partial_{l}A_{m}-B_{l}\partial_{l}A_{i}+\left(A_{l}\partial_{l}\right)B_{i}+\left(B_{l}\partial_{l}\right)A_{i} & = & \,\,\text{(Eq. \ref{EqIndexReplace})}\nonumber \\
\delta_{il}\delta_{jm}A_{j}\partial_{l}B_{m}-\left(A_{l}\partial_{l}\right)B_{i}+\delta_{il}\delta_{jm}B_{j}\partial_{l}A_{m}-\left(B_{l}\partial_{l}\right)A_{i}+\left(A_{l}\partial_{l}\right)B_{i}+\left(B_{l}\partial_{l}\right)A_{i} & = & \,\,\text{(grouping)}\nonumber \\
\delta_{il}\delta_{jm}A_{j}\partial_{l}B_{m}+\delta_{il}\delta_{jm}B_{j}\partial_{l}A_{m} & = & \,\,\text{(cancellation)}\nonumber \\
A_{m}\partial_{i}B_{m}+B_{m}\partial_{i}A_{m} & = & \,\,\text{(Eq. \ref{EqIndexReplace})}\nonumber \\
\partial_{i}\left(A_{m}B_{m}\right) & = & \,\,\text{(product rule)}\nonumber \\
 & = & \left[\nabla\left(\mathbf{A}\cdot\mathbf{B}\right)\right]_{i}\,\text{(Eqs. \ref{eqGrad} \& \ref{eqDiv})}\nonumber
\end{eqnarray}
}Because $i$ is a free index the identity is proved for all components.

$\bullet$ $\nabla\cdot\left(\mathbf{A}\times\mathbf{B}\right)=\mathbf{B}\cdot\left(\nabla\times\mathbf{A}\right)-\mathbf{A}\cdot\left(\nabla\times\mathbf{B}\right)$:
\begin{equation}
\begin{aligned}\nabla\cdot\left(\mathbf{A}\times\mathbf{B}\right) & =\partial_{i}\left[\mathbf{A}\times\mathbf{B}\right]_{i} & \,\,\,\,\,\,\,\,\,\,\,\,\,\, & \text{(Eq. \ref{eqDiv})}\\
 & =\partial_{i}\left(\epsilon_{ijk}A_{j}B_{k}\right) &  & \text{(Eq. \ref{EqCrossProduct})}\\
 & =\epsilon_{ijk}\partial_{i}\left(A_{j}B_{k}\right) &  & \text{(\ensuremath{\partial}\ not acting on \ensuremath{\epsilon})}\\
 & =\epsilon_{ijk}\left(B_{k}\partial_{i}A_{j}+A_{j}\partial_{i}B_{k}\right) &  & \text{(product rule)}\\
 & =\epsilon_{ijk}B_{k}\partial_{i}A_{j}+\epsilon_{ijk}A_{j}\partial_{i}B_{k} &  & \text{(distributivity)}\\
 & =\epsilon_{kij}B_{k}\partial_{i}A_{j}-\epsilon_{jik}A_{j}\partial_{i}B_{k} &  & \text{(Eq. \ref{EqEpsilonCycle})}\\
 & =B_{k}\left(\epsilon_{kij}\partial_{i}A_{j}\right)-A_{j}\left(\epsilon_{jik}\partial_{i}B_{k}\right) &  & \text{(commutativity \& grouping)}\\
 & =B_{k}\left[\nabla\times\mathbf{A}\right]_{k}-A_{j}\left[\nabla\times\mathbf{B}\right]_{j} &  & \text{(Eq. \ref{EqCurl})}\\
 & =\mathbf{B}\cdot\left(\nabla\times\mathbf{A}\right)-\mathbf{A}\cdot\left(\nabla\times\mathbf{B}\right) &  & \text{(Eq. \ref{eqDotProduct})}
\end{aligned}
\end{equation}

$\bullet$ $\nabla\times\left(\mathbf{A}\times\mathbf{B}\right)=\left(\mathbf{B}\cdot\nabla\right)\mathbf{A}+\left(\nabla\cdot\mathbf{B}\right)\mathbf{A}-\left(\nabla\cdot\mathbf{A}\right)\mathbf{B}-\left(\mathbf{A}\cdot\nabla\right)\mathbf{B}$:
\begin{equation}
\begin{aligned}\hspace{-0.5cm}\left[\nabla\times\left(\mathbf{A}\times\mathbf{B}\right)\right]_{i} & =\epsilon_{ijk}\partial_{j}\left[\mathbf{A}\times\mathbf{B}\right]_{k} & \,\,\,\,\, & \text{(Eq. \ref{EqCurl})}\\
 & =\epsilon_{ijk}\partial_{j}\left(\epsilon_{klm}A_{l}B_{m}\right) &  & \text{(Eq. \ref{EqCrossProduct})}\\
 & =\epsilon_{ijk}\epsilon_{klm}\partial_{j}\left(A_{l}B_{m}\right) &  & \text{(\ensuremath{\partial}\ not acting on \ensuremath{\epsilon})}\\
 & =\epsilon_{ijk}\epsilon_{klm}\left(B_{m}\partial_{j}A_{l}+A_{l}\partial_{j}B_{m}\right) &  & \text{(product rule)}\\
 & =\epsilon_{ijk}\epsilon_{lmk}\left(B_{m}\partial_{j}A_{l}+A_{l}\partial_{j}B_{m}\right) &  & \text{(Eq. \ref{EqEpsilonCycle})}\\
 & =\left(\delta_{il}\delta_{jm}-\delta_{im}\delta_{jl}\right)\left(B_{m}\partial_{j}A_{l}+A_{l}\partial_{j}B_{m}\right) &  & \text{(Eq. \ref{EqEpsilonDelta})}\\
 & =\delta_{il}\delta_{jm}B_{m}\partial_{j}A_{l}+\delta_{il}\delta_{jm}A_{l}\partial_{j}B_{m}-\delta_{im}\delta_{jl}B_{m}\partial_{j}A_{l}-\delta_{im}\delta_{jl}A_{l}\partial_{j}B_{m} &  & \text{(distributivity)}\\
 & =B_{m}\partial_{m}A_{i}+A_{i}\partial_{m}B_{m}-B_{i}\partial_{j}A_{j}-A_{j}\partial_{j}B_{i} &  & \text{(Eq. \ref{EqIndexReplace})}\\
 & =\left(B_{m}\partial_{m}\right)A_{i}+\left(\partial_{m}B_{m}\right)A_{i}-\left(\partial_{j}A_{j}\right)B_{i}-\left(A_{j}\partial_{j}\right)B_{i} &  & \text{(grouping)}\\
 & =\left[\left(\mathbf{B}\cdot\nabla\right)\mathbf{A}\right]_{i}+\left[\left(\nabla\cdot\mathbf{B}\right)\mathbf{A}\right]_{i}-\left[\left(\nabla\cdot\mathbf{A}\right)\mathbf{B}\right]_{i}-\left[\left(\mathbf{A}\cdot\nabla\right)\mathbf{B}\right]_{i} &  & \text{(Eqs. \ref{eqANabla} \& \ref{eqDiv})}\\
 & =\left[\left(\mathbf{B}\cdot\nabla\right)\mathbf{A}+\left(\nabla\cdot\mathbf{B}\right)\mathbf{A}-\left(\nabla\cdot\mathbf{A}\right)\mathbf{B}-\left(\mathbf{A}\cdot\nabla\right)\mathbf{B}\right]_{i} &  & \text{(Eq. \ref{eqIndexDistributive1})}
\end{aligned}
\end{equation}
Because $i$ is a free index the identity is proved for all components.

$\bullet$ $\left(\mathbf{A}\times\mathbf{B}\right)\cdot\left(\mathbf{C}\times\mathbf{D}\right)=\begin{vmatrix}\begin{array}{cc}
\mathbf{A}\cdot\mathbf{C} & \mathbf{A}\cdot\mathbf{D}\\
\mathbf{B}\cdot\mathbf{C} & \mathbf{B}\cdot\mathbf{D}
\end{array}\end{vmatrix}$:
\begin{equation}
\begin{aligned}\left(\mathbf{A}\times\mathbf{B}\right)\cdot\left(\mathbf{C}\times\mathbf{D}\right) & =\left[\mathbf{A}\times\mathbf{B}\right]_{i}\left[\mathbf{C}\times\mathbf{D}\right]_{i} & \,\,\,\,\,\, & \text{(Eq. \ref{eqDotProduct})}\\
 & =\epsilon_{ijk}A_{j}B_{k}\epsilon_{ilm}C_{l}D_{m} &  & \text{(Eq. \ref{EqCrossProduct})}\\
 & =\epsilon_{ijk}\epsilon_{ilm}A_{j}B_{k}C_{l}D_{m} &  & \text{(commutativity)}\\
 & =\left(\delta_{jl}\delta_{km}-\delta_{jm}\delta_{kl}\right)A_{j}B_{k}C_{l}D_{m} &  & \text{(Eqs. \ref{EqEpsilonCycle} \& \ref{EqEpsilonDelta})}\\
 & =\delta_{jl}\delta_{km}A_{j}B_{k}C_{l}D_{m}-\delta_{jm}\delta_{kl}A_{j}B_{k}C_{l}D_{m} &  & \text{(distributivity)}\\
 & =\left(\delta_{jl}A_{j}C_{l}\right)\left(\delta_{km}B_{k}D_{m}\right)-\left(\delta_{jm}A_{j}D_{m}\right)\left(\delta_{kl}B_{k}C_{l}\right) &  & \text{(commutativity \& grouping)}\\
 & =\left(A_{l}C_{l}\right)\left(B_{m}D_{m}\right)-\left(A_{m}D_{m}\right)\left(B_{l}C_{l}\right) &  & \text{(Eq. \ref{EqIndexReplace})}\\
 & =\left(\mathbf{A}\cdot\mathbf{C}\right)\left(\mathbf{B}\cdot\mathbf{D}\right)-\left(\mathbf{A}\cdot\mathbf{D}\right)\left(\mathbf{B}\cdot\mathbf{C}\right) &  & \text{(Eq. \ref{eqDotProduct})}\\
 & =\begin{vmatrix}\begin{array}{cc}
\mathbf{A}\cdot\mathbf{C} & \mathbf{A}\cdot\mathbf{D}\\
\mathbf{B}\cdot\mathbf{C} & \mathbf{B}\cdot\mathbf{D}
\end{array}\end{vmatrix} &  & \text{(definition of determinant)}
\end{aligned}
\end{equation}

$\bullet$ $\left(\mathbf{A}\times\mathbf{B}\right)\times\left(\mathbf{C}\times\mathbf{D}\right)=\left[\mathbf{D}\cdot\left(\mathbf{A}\times\mathbf{B}\right)\right]\mathbf{C}-\left[\mathbf{C}\cdot\left(\mathbf{A}\times\mathbf{B}\right)\right]\mathbf{D}$:
\begin{equation}
\begin{aligned}\left[\left(\mathbf{A}\times\mathbf{B}\right)\times\left(\mathbf{C}\times\mathbf{D}\right)\right]_{i} & =\epsilon_{ijk}\left[\mathbf{A}\times\mathbf{B}\right]_{j}\left[\mathbf{C}\times\mathbf{D}\right]_{k} & \,\,\,\,\, & \text{(Eq. \ref{EqCrossProduct})}\\
 & =\epsilon_{ijk}\epsilon_{jmn}A_{m}B_{n}\epsilon_{kpq}C_{p}D_{q} &  & \text{(Eq. \ref{EqCrossProduct})}\\
 & =\epsilon_{ijk}\epsilon_{kpq}\epsilon_{jmn}A_{m}B_{n}C_{p}D_{q} &  & \text{(commutativity)}\\
 & =\epsilon_{ijk}\epsilon_{pqk}\epsilon_{jmn}A_{m}B_{n}C_{p}D_{q} &  & \text{(Eq. \ref{EqEpsilonCycle})}\\
 & =\left(\delta_{ip}\delta_{jq}-\delta_{iq}\delta_{jp}\right)\epsilon_{jmn}A_{m}B_{n}C_{p}D_{q} &  & \text{(Eq. \ref{EqEpsilonDelta})}\\
 & =\left(\delta_{ip}\delta_{jq}\epsilon_{jmn}-\delta_{iq}\delta_{jp}\epsilon_{jmn}\right)A_{m}B_{n}C_{p}D_{q} &  & \text{(distributivity)}\\
 & =\left(\delta_{ip}\epsilon_{qmn}-\delta_{iq}\epsilon_{pmn}\right)A_{m}B_{n}C_{p}D_{q} &  & \text{(Eq. \ref{EqIndexReplace})}\\
 & =\delta_{ip}\epsilon_{qmn}A_{m}B_{n}C_{p}D_{q}-\delta_{iq}\epsilon_{pmn}A_{m}B_{n}C_{p}D_{q} &  & \text{(distributivity)}\\
 & =\epsilon_{qmn}A_{m}B_{n}C_{i}D_{q}-\epsilon_{pmn}A_{m}B_{n}C_{p}D_{i} &  & \text{(Eq. \ref{EqIndexReplace})}\\
 & =\epsilon_{qmn}D_{q}A_{m}B_{n}C_{i}-\epsilon_{pmn}C_{p}A_{m}B_{n}D_{i} &  & \text{(commutativity)}\\
 & =\left(\epsilon_{qmn}D_{q}A_{m}B_{n}\right)C_{i}-\left(\epsilon_{pmn}C_{p}A_{m}B_{n}\right)D_{i} &  & \text{(grouping)}\\
 & =\left[\mathbf{D}\cdot\left(\mathbf{A}\times\mathbf{B}\right)\right]C_{i}-\left[\mathbf{C}\cdot\left(\mathbf{A}\times\mathbf{B}\right)\right]D_{i} &  & \text{(Eq. \ref{EqScalarTripleProduct})}\\
 & =\left[\left[\mathbf{D}\cdot\left(\mathbf{A}\times\mathbf{B}\right)\right]\mathbf{C}\right]_{i}-\left[\left[\mathbf{C}\cdot\left(\mathbf{A}\times\mathbf{B}\right)\right]\mathbf{D}\right]_{i} &  & \text{(definition of index)}\\
 & =\left[\left[\mathbf{D}\cdot\left(\mathbf{A}\times\mathbf{B}\right)\right]\mathbf{C}-\left[\mathbf{C}\cdot\left(\mathbf{A}\times\mathbf{B}\right)\right]\mathbf{D}\right]_{i} &  & \text{(Eq. \ref{eqIndexDistributive1})}
\end{aligned}
\end{equation}
Because $i$ is a free index the identity is proved for all components.

\pagebreak{}

\section{Metric Tensor\label{secMetricTensor}}

$\bullet$ This is a rank-2 tensor which may also be called the fundamental
tensor.

$\bullet$ The main purpose of the metric tensor is to generalize
the concept of distance to general curvilinear coordinate frames and
maintain the invariance of distance in different coordinate systems.

$\bullet$ In orthonormal Cartesian coordinate systems the distance
element squared, $\left(ds\right)^{2}$, between two infinitesimally
neighboring points in space, one with coordinates $x^{i}$ and the
other with coordinates $x^{i}+dx^{i}$, is given by
\begin{equation}
\left(ds\right)^{2}=dx^{i}dx^{i}=\delta_{ij}dx^{i}dx^{j}
\end{equation}
This definition of distance is the key to introducing a rank-2 tensor,
$g_{ij}$, called the metric tensor which, for a general coordinate
system, is defined by
\begin{equation}
\left(ds\right)^{2}=g_{ij}dx^{i}dx^{j}
\end{equation}
The metric tensor has also a contravariant form, i.e. $g^{ij}$.

$\bullet$ The components of the metric tensor are given by:
\begin{equation}
g_{ij}=\mathbf{E}_{i}\cdot\mathbf{E}_{j}\,\,\,\,\,\,\,\,\,\,\,\,\,\,\,\,\,\,\,\,\&\,\,\,\,\,\,\,\,\,\,\,\,\,\,\,\,\,\,\,\,g^{ij}=\mathbf{E}^{i}\cdot\mathbf{E}^{j}
\end{equation}

where the indexed $\mathbf{E}$ are the covariant and contravariant
basis vectors as defined in $\S$ \ref{subCovariantContravariant}.

$\bullet$ The mixed type metric tensor is given by:
\begin{equation}
g_{\,\,j}^{i}=\mathbf{E}^{i}\cdot\mathbf{E}_{j}=\delta_{\,\,j}^{i}\,\,\,\,\,\,\,\,\,\,\,\,\,\,\,\,\,\,\,\,\&\,\,\,\,\,\,\,\,\,\,\,\,\,\,\,\,\,\,\,\,g_{i}^{\,\,j}=\mathbf{E}_{i}\cdot\mathbf{E}^{j}=\delta_{i}^{\,\,j}
\end{equation}
and hence it is the same as the unity tensor.

$\bullet$ For a coordinate system in which the metric tensor can
be cast in a diagonal form where the diagonal elements are $\pm1$
the metric is called flat.

$\bullet$ For Cartesian coordinate systems, which are orthonormal
flat-space systems, we have
\begin{equation}
g^{ij}=\delta^{ij}=g_{ij}=\delta_{ij}
\end{equation}

$\bullet$ The metric tensor is symmetric, that is
\begin{equation}
g_{ij}=g_{ji}\,\,\,\,\,\,\,\,\,\,\,\,\,\&\,\,\,\,\,\,\,\,\,\,\,g^{ij}=g^{ji}
\end{equation}

$\bullet$ The contravariant metric tensor is used for raising indices
of covariant tensors and the covariant metric tensor is used for lowering
indices of contravariant tensors, e.g.
\begin{equation}
A^{i}=g^{ij}A_{j}\,\,\,\,\,\,\,\,\,\,\,\,\,\,\,\,\,\,\,\,A_{i}=g_{ij}A^{j}
\end{equation}
where the metric tensor acts, like a Kronecker delta, as an index
replacement operator. Hence, any tensor can be cast into a covariant
or a contravariant form, as well as a mixed form. \textcolor{black}{However,
the order of the indices should be respected in this process, e.g.}
\begin{equation}
A_{\,\,j}^{i}=g_{jk}A^{ik}\ne A_{j}^{\,\,\,i}=g_{jk}A^{ki}
\end{equation}
Some authors insert dots (e.g. $A_{j}^{\cdot\,i}$) to remove any
ambiguity about the order of the indices.

$\bullet$ The covariant and contravariant metric tensors are inverses
of each other, that is
\begin{equation}
\left[g_{ij}\right]=\left[g^{ij}\right]^{-1}\,\,\,\,\,\,\,\,\,\,\&\,\,\,\,\,\,\,\,\,\,\left[g^{ij}\right]=\left[g_{ij}\right]^{-1}
\end{equation}
Hence
\begin{equation}
g^{ik}g_{kj}=\delta_{\,\,j}^{i}\,\,\,\,\,\,\,\,\,\,\&\,\,\,\,\,\,\,\,\,\,g_{ik}g^{kj}=\delta_{i}^{\,\,j}
\end{equation}

$\bullet$ It is common to reserve the ``metric tensor'' to the
covariant form and call the contravariant form, which is its inverse,
the ``associate'' or ``conjugate'' or ``reciprocal'' metric
tensor.

$\bullet$ As a tensor, the metric has a significance regardless of
any coordinate system although it requires a coordinate system to
be represented in a specific form.

$\bullet$ For orthogonal coordinate systems the metric tensor is
diagonal, i.e. $g_{ij}=g^{ij}=0$ for $i\ne j$.

$\bullet$ For flat-space orthonormal Cartesian coordinate systems
in a 3D space, the metric tensor is given by:
\begin{equation}
\left[g_{ij}\right]=\left[\delta_{ij}\right]=\left[\begin{array}{ccc}
1 & 0 & 0\\
0 & 1 & 0\\
0 & 0 & 1
\end{array}\right]=\left[\delta^{ij}\right]=\left[g^{ij}\right]
\end{equation}

$\bullet$ For cylindrical coordinate systems with coordinates ($\rho,\phi,z$),
the metric tensor is given by:
\begin{equation}
\left[g_{ij}\right]=\left[\begin{array}{ccc}
1 & 0 & 0\\
0 & \rho^{2} & 0\\
0 & 0 & 1
\end{array}\right]\,\,\,\,\,\,\,\,\,\,\&\,\,\,\,\,\,\,\,\,\,\left[g^{ij}\right]=\left[\begin{array}{ccc}
1 & 0 & 0\\
0 & \frac{1}{\rho^{2}} & 0\\
0 & 0 & 1
\end{array}\right]
\end{equation}

$\bullet$ For spherical coordinate systems with coordinates ($r,\theta,\phi$),
the metric tensor is given by:
\begin{equation}
\left[g_{ij}\right]=\left[\begin{array}{ccc}
1 & 0 & 0\\
0 & r^{2} & 0\\
0 & 0 & r^{2}\sin^{2}\theta
\end{array}\right]\,\,\,\,\,\,\,\,\,\,\&\,\,\,\,\,\,\,\,\,\,\left[g^{ij}\right]=\left[\begin{array}{ccc}
1 & 0 & 0\\
0 & \frac{1}{r^{2}} & 0\\
0 & 0 & \frac{1}{r^{2}\sin^{2}\theta}
\end{array}\right]
\end{equation}

\pagebreak{}

\section{Covariant Differentiation\label{secCovariantDifferentiation}}

$\bullet$ The ordinary derivative of a tensor is not a tensor in
general. The objective of covariant differentiation is to ensure the
invariance of derivative (i.e. being a tensor) in general coordinate
systems, and this results in applying more sophisticated rules using
Christoffel symbols where different differentiation rules for covariant
and contravariant indices apply. The resulting covariant derivative
is a tensor which is one rank higher than the differentiated tensor.

$\bullet$ Christoffel symbol of the second kind is defined by:
\begin{equation}
\left\{ _{ij}^{k}\right\} =\frac{g^{kl}}{2}\left(\frac{\partial g_{il}}{\partial x^{j}}+\frac{\partial g_{jl}}{\partial x^{i}}-\frac{\partial g_{ij}}{\partial x^{l}}\right)
\end{equation}
where the indexed $g$ is the metric tensor in its contravariant and
covariant forms with implied summation over $l$. It is noteworthy
that Christoffel symbols are not tensors.

$\bullet$ The Christoffel symbols of the second kind are symmetric
in their two lower indices:
\begin{equation}
\left\{ _{ij}^{k}\right\} =\left\{ _{ji}^{k}\right\}
\end{equation}

$\bullet$ For Cartesian coordinate systems, the Christoffel symbols
are zero for all the values of indices.

$\bullet$ For cylindrical coordinate systems ($\rho,\phi,z$), the
Christoffel symbols are zero for all the values of indices except:
\begin{eqnarray}
\left\{ _{22}^{1}\right\}  & = & -\rho\\
\left\{ _{12}^{2}\right\}  & = & \begin{aligned}\left\{ _{21}^{2}\right\}  & \begin{aligned}\,\,=\,\,\,\, & \frac{1}{\rho}\end{aligned}
\end{aligned}
\nonumber
\end{eqnarray}
where ($1,2,3$) stand for ($\rho,\phi,z$).

$\bullet$ For spherical coordinate systems ($r,\theta,\phi$), the
Christoffel symbols are zero for all the values of indices except:
\begin{eqnarray}
\left\{ _{22}^{1}\right\}  & = & -r\\
\left\{ _{33}^{1}\right\}  & = & -r\sin^{2}\theta\nonumber \\
\left\{ _{12}^{2}\right\}  & = & \begin{aligned}\left\{ _{21}^{2}\right\}  & \begin{aligned}\,=\,\,\,\,\,\, & \frac{1}{r}\end{aligned}
\end{aligned}
\nonumber \\
\left\{ _{33}^{2}\right\}  & = & -\sin\theta\cos\theta\nonumber \\
\left\{ _{13}^{3}\right\}  & = & \begin{aligned}\begin{aligned}\left\{ _{31}^{3}\right\}  & =\end{aligned}
 & \frac{1}{r}\end{aligned}
\nonumber \\
\left\{ _{23}^{3}\right\}  & = & \begin{aligned}\begin{aligned}\left\{ _{32}^{3}\right\}  & =\end{aligned}
 & \cot\theta\end{aligned}
\nonumber
\end{eqnarray}
where ($1,2,3$) stand for ($r,\theta,\phi$).

$\bullet$ For a differentiable scalar $f$ the covariant derivative
is the same as the normal partial derivative, that is:
\begin{equation}
f_{;i}=f_{,i}=\partial_{i}f
\end{equation}
This is justified by the fact that the covariant derivative is different
from the normal partial derivative because the basis vectors in general
coordinate systems are dependent on their spatial position, and since
a scalar is independent of the basis vectors the covariant and partial
derivatives are identical.

$\bullet$ For a differentiable vector $\mathbf{A}$ the covariant
derivative is:
\begin{equation}
\begin{aligned}A_{j;i} & =\partial_{i}A_{j}-\left\{ _{ji}^{k}\right\} A_{k} & \,\,\,\,\,\,\,\,\,\,\,\,\,\,\, & \text{(covariant)}\\
A_{\,\,;i}^{j} & =\partial_{i}A^{j}+\left\{ _{ki}^{j}\right\} A^{k} &  & \text{(contravariant)}
\end{aligned}
\end{equation}

$\bullet$ For a differentiable rank-2 tensor $\mathbf{A}$ the covariant
derivative is:
\begin{equation}
\begin{aligned}A_{jk;i} & =\partial_{i}A_{jk}-\left\{ _{ji}^{l}\right\} A_{lk}-\left\{ _{ki}^{l}\right\} A_{jl} & \,\,\,\,\,\,\,\,\,\,\,\,\, & \text{(covariant)}\\
A_{\,\,\,;i}^{jk} & =\partial_{i}A^{jk}+\left\{ _{li}^{j}\right\} A^{lk}+\left\{ _{li}^{k}\right\} A^{jl} &  & \text{(contravariant)}\\
A_{j;i}^{k} & =\partial_{i}A_{j}^{k}+\left\{ _{li}^{k}\right\} A_{j}^{l}-\left\{ _{ji}^{l}\right\} A_{l}^{k} &  & \text{(mixed)}
\end{aligned}
\end{equation}

$\bullet$ For a differentiable rank-$n$ tensor $\mathbf{A}$ the
covariant derivative is:
\begin{eqnarray}
A_{lm\ldots p;q}^{ij\ldots k} & = & \partial_{q}A_{lm\ldots p}^{ij\ldots k}+\left\{ _{aq}^{i}\right\} A_{lm\ldots p}^{aj\ldots k}+\left\{ _{aq}^{j}\right\} A_{lm\ldots p}^{ia\ldots k}+\cdots+\left\{ _{aq}^{k}\right\} A_{lm\ldots p}^{ij\ldots a}\\
 &  & \,\,\,\,\,\,\,\,\,\,\,\,\,\,\,\,\,\,\,\,\,\,-\left\{ _{lq}^{a}\right\} A_{am\ldots p}^{ij\ldots k}-\left\{ _{mq}^{a}\right\} A_{la\ldots p}^{ij\ldots k}-\cdots-\left\{ _{pq}^{a}\right\} A_{lm\ldots a}^{ij\ldots k}\nonumber
\end{eqnarray}

$\bullet$ From the last three points a pattern for covariant differentiation
emerges, that is it starts with a partial derivative term then for
each tensor index an extra Christoffel symbol term is added, positive
for superscripts and negative for subscripts, where the differentiation
index is the second of the lower indices in the Christoffel symbol.

$\bullet$ Since the Christoffel symbols are identically zero in Cartesian
coordinate systems, the covariant derivative is the same as the normal
partial derivative for all tensor ranks.

$\bullet$ The covariant derivative of the metric tensor is zero in
all coordinate systems.

$\bullet$ Several rules of normal differentiation similarly apply
to covariant differentiation. For example, covariant differentiation
is a linear operation with respect to algebraic sums of tensor terms:
\begin{equation}
\partial_{;i}\left(a\mathbf{A}\pm b\mathbf{B}\right)=a\partial_{;i}\mathbf{A}\pm b\partial_{;i}\mathbf{B}
\end{equation}
where $a$ and $b$ are scalar constants and $\mathbf{A}$ and $\mathbf{B}$
are differentiable tensor fields. The product rule of normal differentiation
also applies to covariant differentiation of tensor multiplication:
\begin{equation}
\partial_{;i}\left(\mathbf{A}\mathbf{B}\right)=\left(\partial_{;i}\mathbf{A}\right)\mathbf{B}+\mathbf{A}\partial_{;i}\mathbf{B}
\end{equation}
This rule is also valid for the inner product of tensors because the
inner product is an outer product operation followed by a contraction
of indices, and covariant differentiation and contraction of indices
commute.

$\bullet$ The covariant derivative operator can bypass the raising/lowering
index operator:
\begin{equation}
A_{i}=g_{ij}A^{j}\,\,\,\,\,\,\,\,\,\,\,\,\,\,\,\Longrightarrow\,\,\,\,\,\,\,\,\,\,\,\,\,\,\,\partial_{;m}A_{i}=g_{ij}\partial_{;m}A^{j}
\end{equation}
and hence the metric behaves like a constant with respect to the covariant
operator.

$\bullet$ A principal difference between normal partial differentiation
and covariant differentiation is that for successive differential
operations the partial derivative operators do commute with each other
(assuming certain continuity conditions) but the covariant operators
do not commute, that is
\begin{equation}
\partial_{i}\partial_{j}=\partial_{j}\partial_{i}\,\,\,\,\,\,\,\,\,\,\,\,\,\,\,\text{but \,\,\,\,\,\,\,\,\,\,\,\,\,\,\,}\partial_{;i}\partial_{;j}\ne\partial_{;j}\partial_{;i}
\end{equation}

$\bullet$ Higher order covariant derivatives are similarly defined
as derivatives of derivatives; however the order of differentiation
should be respected (refer to the previous point).

\pagebreak{}

\phantomsection
\addcontentsline{toc}{section}{References}
\bibliographystyle{unsrt}

Note: As well as the references cited above, I benefited during the
writing of these notes from many sources such as tutorials, presentations,
and articles which I found on the Internet authored or composed by
other people. As it is difficult or impossible to retrace and state
all these sources, I make a general acknowledgment to all those who
made their documents available to the public.
\end{document}